\documentclass[graybox]{svmult}

\usepackage{type1cm}        
\usepackage{makeidx}         
\usepackage{graphicx}        
\usepackage{multicol}        
\usepackage[bottom]{footmisc}

\usepackage{newtxtext}       %
\usepackage[varvw]{newtxmath}       

\usepackage[x11names]{xcolor}
\definecolor{colorLevel1}{rgb}{0.1804,0.4353,0.5569}
\definecolor{colorLevel3}{rgb}{0.1608,0.6863,0.4980}
\definecolor{colorLevel2}{rgb}{0.7412,0.8745,0.1490}
\colorlet{colorShading}{gray} 
\colorlet{colorShadingLight}{lightgray} 
\colorlet{colorMesh}{black!80} 
\colorlet{colorTree}{red}  
\colorlet{colorDirichlet}{blue} 
\colorlet{colorInterface}{blue} 
\colorlet{colorDeactivated}{cyan} 
\colorlet{colorHomology}{orange} 

\usepackage{tikz}
\usepackage{pgfplots}
\usetikzlibrary{external}
\tikzexternaldisable

\usepackage{lipsum}
\usepackage{amsfonts}
\usepackage{graphicx}
\usepackage{epstopdf}
\usepackage{algorithmic}
\ifpdf
  \DeclareGraphicsExtensions{.eps,.pdf,.png,.jpg}
\else
  \DeclareGraphicsExtensions{.eps}
\fi
\usepackage{subcaption}
\usepackage{amscd}
\usepackage{cleveref}
\usepackage{algorithm}

\newtheorem{assumption}{Assumption}

\title{Tree-cotree gauging for hierarchical B-splines\thanks{Submitted to the editors DATE.
\funding{TODO}}}

\author{Melina Merkel\thanks{Computational Electromagnetics Group and
              Centre for Computational Engineering,
              Technische Universität Darmstadt,
              Darmstadt, Germany,
              (\email{melina.merkel@tu-darmstadt.de}).}
\and Rafael Vázquez\thanks{Departamento de Matem\'atica Aplicada,
              Universidade de Santiago de Compostela,
              Santiago de Compostela, Spain
              and
              Galician Centre for Mathematical Research and Technology (CITMAga),
              Santiago de Compostela, Spain, (\email{rafael.vazquez@usc.es}).}
}

\usepackage{siunitx}
\usepackage[
    nogroupskip,
    acronym,
    nomain,
    xindy,
    stylemods={longbooktabs},
    indexonlyfirst,
    sort=use,
    automake
]
{glossaries-extra}
\setglossarystyle{index} 

\usepackage{amsopn}

\newcommand{\supp}{\mathrm{supp}\,}
\newcommand{\map}{\mathbf{F}}
\newcommand{\grad}{\mathbf{grad}\,}
\newcommand{\graph}{\mathcal{G}}
\newcommand{\tree}{\mathcal{T}}

\newcommand{\nodeSet}{\mathcal{N}}
\newcommand{\edgeSet}{\mathcal{E}}
\newcommand{\basis}{\mathcal{B}}
\newcommand{\hbasis}{\mathcal{H}}
\let\hat\widehat

\newcommand{\mm}[1]{\textcolor{black}{#1}}
\newcommand{\rv}[1]{\textcolor{black}{#1}}
\newcommand{\graphl}{\ensuremath{\graph_{\ell}}}
\newcommand{\act}{\ensuremath{\mathcal{A}}}
\newcommand{\deactivated}{\ensuremath{\mathcal{D}}}

\newcommand{\edgesboundary}[1][\ell]{\ensuremath{\edgeSet_{#1}^{\partial}}}
\newcommand{\edgesinterior}[1][\ell]{\ensuremath{\edgeSet_{#1}^{\mathrm{int}}}}
\newcommand{\edgesactive}[1][\ell]{\ensuremath{\edgeSet_{#1}^{\act}}}
\newcommand{\edgesdeactivated}[1][\ell]{\ensuremath{\edgeSet_{#1}^{\deactivated}}}

\newcommand{\nodesboundary}[1][\ell]{\ensuremath{\nodeSet_{#1}^{\partial}}}
\newcommand{\nodesinterior}[1][\ell]{\ensuremath{\nodeSet_{#1}^{\mathrm{int}}}}
\newcommand{\nodesactive}[1][\ell]{\ensuremath{\nodeSet_{#1}^{\act}}}
\newcommand{\nodesdeactivated}[1][\ell]{\ensuremath{\nodeSet_{#1}^{\deactivated}}}

\newcommand{\treeactive}[1][\ell]{\ensuremath{\tree_{#1}^{\act}}}
\newcommand{\treedeactivated}[1][\ell]{\ensuremath{\tree_{#1}^{\deactivated}}}
\newcommand{\multileveltree}{\ensuremath{\mathcal{M}}}

\newcommand{\curl}{\mathbf{curl}}

\newcommand{\Hcurl}      {\ensuremath{\mathbf{H}(\mathbf{curl};\Omega)}} 

\newcommand{\field}[1]{\ensuremath{\mathbf{#1}}}
\newcommand{\Afield}{\ensuremath{\field{A}}}

\newcommand{\Efield}{\ensuremath{\field{E}}}
\newcommand{\Jsrc}{\ensuremath{\field{J}_{\mathrm{src}}}}

\GlsXtrEnablePreLocationTag{page~}{pages~}

\newglossary[slg]{symbolslist}{syi}{syg}{List of Symbols}

\makeglossaries

\newglossaryentry{knotvec}{
    name={\ensuremath{\Xi}},
    description={Knot vector},
    type=symbolslist,
    sort={},
}
\newglossaryentry{Omegahat}{
    name={\ensuremath{\hat \Omega}},
    description={Parametric domain},
    type=symbolslist,
    sort={},
}
\newglossaryentry{Sp}{
    name={\ensuremath{S_p(\Xi)}},
    description={Univariate B-spline basis of degree $p$},
    type=symbolslist,
    sort={},
}
\newglossaryentry{Xkhat}{
    name={\ensuremath{\hat{X}^k}},
    description={B-spline spaces of de Rham sequence for $k = 0,\ldots, 3$ in the parametric domain},
    type=symbolslist,
    sort={},
}
\newglossaryentry{Bkhat}{
    name={\ensuremath{\hat \basis^k}},
    description={The basis for $\hat{X}^k$, for $k = 0,\ldots, 3$ in the parametric domain},
    type=symbolslist,
    sort={},
}
\newglossaryentry{Omega}{
    name={\ensuremath{\Omega}},
    description={Physical domain},
    type=symbolslist,
    sort={},
}
\newglossaryentry{F}{
    name={\ensuremath{\map}},
    description={Parametrization of the physical domain},
    type=symbolslist,
    sort={},
}
\newglossaryentry{Xk}{
    name={\ensuremath{X^k}},
    description={B-spline spaces of de Rham sequence for $k = 0,\ldots, 3$ in the physical domain},
    type=symbolslist,
    sort={},
}
\newglossaryentry{DF}{
    name={\ensuremath{D\map}},
    description={Jacobian of the parametrization $\gls{F}$},
    type=symbolslist,
    sort={},
}
\newglossaryentry{Bk}{
    name={\ensuremath{\basis^k}},
    description={The basis for ${X}^k$, for $k = 0,\ldots, 3$ in the physical domain},
    type=symbolslist,
    sort={},
}
\newglossaryentry{Omegai}{
    name={\ensuremath{\Omega^{(i)}}},
    description={Patch $i$},
    type=symbolslist,
    sort={},
}
\newglossaryentry{Fi}{
    name={\ensuremath{\map^{(i)}}},
    description={Parametrization of the $i$-th patch},
    type=symbolslist,
    sort={},
}
\newglossaryentry{Pi}{
    name={\ensuremath{P_i}},
    description={$i$-th Greville point},
    type=symbolslist,
    sort={},
}
\newglossaryentry{G}{
    name={\ensuremath{G}},
    description={Greville grid},
    type=symbolslist,
    sort={},
}
\newglossaryentry{Zk}{
    name={\ensuremath{Z^k}},
    description={Lowest order finite element spaces of de Rham sequence for $k = 0,\ldots, 3$},
    type=symbolslist,
    sort={},
}
\newglossaryentry{Ik}{
    name={\ensuremath{I^K}},
    description={Isomorphisms between B-spline \gls{Xk} and \gls{Zk}},
    type=symbolslist,
    sort={},
}
\newglossaryentry{l}{
    name={\ensuremath{\ell}},
    description={Refinement level of hierarchical B-splines},
    type=symbolslist,
    sort={},
}
\newglossaryentry{Xkl}{
    name={\ensuremath{X^k_\ell}},
    description={B-spline spaces of de Rham sequence of level $\ell = 0, \ldots, L$ for $k = 0,\ldots, 3$},
    type=symbolslist,
    sort={},
}
\newglossaryentry{Bkl}{
    name={\ensuremath{\basis_\ell^k}},
    description={The basis for ${X}^k_\ell$ of level $\ell$},
    type=symbolslist,
    sort={},
}
\newglossaryentry{Gl}{
    name={\ensuremath{G_\ell}},
    description={Greville grid of level $\ell$},
    type=symbolslist,
    sort={},
}
\newglossaryentry{Zkl}{
    name={\ensuremath{Z^k_\ell}},
    description={Lowest order finite element spaces defined on $G_\ell$ for $k = 0,\ldots, 3$},
    type=symbolslist,
    sort={},
}
\newglossaryentry{Omegal}{
    name={\ensuremath{\Omega_\ell}},
    description={Closed subdomain of level $\ell$},
    type=symbolslist,
    sort={},
}
\newglossaryentry{Dl}{
    name={\ensuremath{D_\ell}},
    description={Region $D_\ell = \Omega_\ell \setminus \Omega_{\ell+1}$},
    type=symbolslist,
    sort={},
}
\newglossaryentry{lp}{
    name={\ensuremath{\ell'}},
    description={Refinement level $\ell' \in \left\{\ell, \ell+1\right\}$},
    type=symbolslist,
    sort={},
}
\newglossaryentry{Bkll}{
    name={\ensuremath{\basis^k_{\ell, \ell'}}},
    description={Subset of basis functions of $\gls{Bkl}$ with support completely contained in $\Omega_{\ell'}$},
    type=symbolslist,
    sort={},
}
\newglossaryentry{HkL}{
    name={\ensuremath{\hbasis^k_{L}}},
    description={Hierarchical B-spline basis of level $L$ for $k = 0,\ldots, 3$},
    type=symbolslist,
    sort={},
}
\newglossaryentry{Hkl}{
    name={\ensuremath{\hbasis^k_{\ell}}},
    description={Hierarchical B-spline basis up to level $\ell$ for $k = 0,\ldots, 3$},
    type=symbolslist,
    sort={},
}
\newglossaryentry{Akl}{
    name={\ensuremath{\act^k_{\ell}}},
    description={Set of active basis functions on level $\ell$},
    type=symbolslist,
    sort={},
}
\newglossaryentry{Dkl}{
    name={\ensuremath{\deactivated^k_{\ell}}},
    description={Set of deactivated basis functions on level $\ell$},
    type=symbolslist,
    sort={},
}
\newglossaryentry{Wkl}{
    name={\ensuremath{W^k_\ell}},
    description={Space of hierarchical B-splines up to level $\ell$, $W^k_\ell = \mathrm{span} \{\hbasis^k_{\ell}\}$},
    type=symbolslist,
    sort={},
}
\newglossaryentry{Xkll}{
    name={\ensuremath{X^k_{\ell,\ell'}}},
    description={Spaces of hierarchical B-splines $\basis^k_{\ell}$ with support completely contained in $\Omega_{\ell'}$, $X^k_{\ell,\ell'} = \mathrm{span} \{\basis^k_{\ell,\ell'}\}$},
    type=symbolslist,
    sort={},
}
\newglossaryentry{Gll}{
    name={\ensuremath{G_{\ell,\ell'}}},
    description={Greville subgrid formed by cells associated to $\basis^3_{\ell, \ell'}$, also used for domain covered by $G_{\ell,\ell'}$},
    type=symbolslist,
    sort={},
}
\newglossaryentry{dGll}{
    name={\ensuremath{\partial G_{\ell,\ell'}}},
    description={Boundary of $\partial G_{\ell,\ell'}$},
    type=symbolslist,
    sort={},
}
\newglossaryentry{DGl}{
    name={\ensuremath{D^G_\ell}},
    description={Domain covered by Greville subgrid $G_{\ell,\ell} \setminus G_{\ell,\ell+1}$ of level $\ell$},
    type=symbolslist,
    sort={},
}
\newglossaryentry{Zkll}{
    name={\ensuremath{Z^k_{\ell,\ell'}}},
    description={Lowest order finite element spaces on $G_{\ell,\ell'}$},
    type=symbolslist,
    sort={},
}
\newglossaryentry{Hk}{
    name={\ensuremath{H^k}},
    description={Cohomology group \mm{@Rafa?}},
    type=symbolslist,
    sort={},
}
\newglossaryentry{graph}{
    name={\ensuremath{\graph}},
    description={Graph from nodes and edges of the Greville grid $G$},
    type=symbolslist,
    sort={},
}
\newglossaryentry{N}{
    name={\ensuremath{\nodeSet}},
    description={Set of nodes in the graph $\graph$},
    type=symbolslist,
    sort={},
}
\newglossaryentry{E}{
    name={\ensuremath{\edgeSet}},
    description={Set of edges in the graph $\graph$},
    type=symbolslist,
    sort={},
}
\newglossaryentry{tree}{
    name={\ensuremath{\tree}},
    description={Spanning tree constructed on the graph $\graph$},
    type=symbolslist,
    sort={},
}
\newglossaryentry{graphl}{
    name={\ensuremath{\graphl}},
    description={Graph from nodes and edges of the Greville subgrid $G_{\ell,\ell}$},
    type=symbolslist,
    sort={},
}
\newglossaryentry{Nl}{
    name={\ensuremath{\nodeSet_{\ell}}},
    description={Set of nodes in the graph $\graphl$},
    type=symbolslist,
    sort={},
}
\newglossaryentry{El}{
    name={\ensuremath{\edgeSet_{\ell}}},
    description={Set of edges in the graph $\graphl$},
    type=symbolslist,
    sort={},
}
\newglossaryentry{Nintl}{
    name={\ensuremath{\nodesinterior}},
    description={Nodes in the interior of $G_{\ell,\ell}$},
    type=symbolslist,
    sort={},
}
\newglossaryentry{nNintl}{
    name={\ensuremath{\# \nodesinterior}},
    description={Number of nodes in the set$\nodesinterior$},
    type=symbolslist,
    sort={},
}
\newglossaryentry{Ndl}{
    name={\ensuremath{\nodesboundary}},
    description={Nodes on the boundary $\partial G_{\ell,\ell}$},
    type=symbolslist,
    sort={},
}
\newglossaryentry{NAl}{
    name={\ensuremath{\nodesactive}},
    description={Interior nodes of $G_{\ell,\ell}$ associated to active basis functions $\act^0_{\ell}$},
    type=symbolslist,
    sort={},
}
\newglossaryentry{NDl}{
    name={\ensuremath{\nodesdeactivated}},
    description={Interior nodes of $G_{\ell,\ell}$ associated to deactivated basis functions $\deactivated^0_{\ell}$},
    type=symbolslist,
    sort={},
}
\newglossaryentry{Eintl}{
    name={\ensuremath{\edgesinterior}},
    description={Edges in the interior of $G_{\ell,\ell}$},
    type=symbolslist,
    sort={},
}
\newglossaryentry{Edl}{
    name={\ensuremath{\edgesboundary}},
    description={Edges on the boundary $\partial G_{\ell,\ell}$},
    type=symbolslist,
    sort={},
}
\newglossaryentry{EAl}{
    name={\ensuremath{\edgesactive}},
    description={Interior edges of $G_{\ell,\ell}$ associated to active basis functions $\act^1_{\ell}$},
    type=symbolslist,
    sort={},
}
\newglossaryentry{EDl}{
    name={\ensuremath{\edgesdeactivated}},
    description={Interior edges of $G_{\ell,\ell}$ associated to deactivated basis functions $\deactivated^1_{\ell}$},
    type=symbolslist,
    sort={},
}
\newglossaryentry{graphdl}{
    name={\ensuremath{\graphl^\partial}},
    description={Subgraph constructed from from $\graphl^\partial = (\nodesboundary,\edgesboundary)$},
    type=symbolslist,
    sort={},
}
\newglossaryentry{Tl}{
    name={\ensuremath{\tree_{\ell}}},
    description={Spanning tree of level $\ell$ constructed on the graph $\graphl$},
    type=symbolslist,
    sort={},
}
\newglossaryentry{nTl}{
    name={},
    description={Number of interior edges in $\tree_{\ell}$},
    type=symbolslist,
    sort={},
}
\newglossaryentry{TAl}{
    name={\ensuremath{\treeactive}},
    description={Restriction of $\tree_{\ell}$ to edges associated to active functions of level $\ell$},
    type=symbolslist,
    sort={},
}
\newglossaryentry{TDl}{
    name={\ensuremath{\treedeactivated}},
    description={Restriction of $\tree_{\ell}$ to edges associated to deactivated functions of level $\ell$},
    type=symbolslist,
    sort={},
}
\newglossaryentry{ML}{
    name={\ensuremath{\multileveltree_{L}}},
    description={Multi-level tree on level $L$},
    type=symbolslist,
    sort={},
}
\newglossaryentry{Ml}{
    name={\ensuremath{\multileveltree_{\ell}}},
    description={Multi-level tree up to level $\ell$},
    type=symbolslist,
    sort={},
}
\makeindex             

\begin{document}

\title*{Tree-cotree gauging for two-dimensional hierarchical splines}
\author{M. Merkel \orcidID{0000-0002-2104-9167} and R. Vázquez}
\institute{M. Merkel \at Computational Electromagnetics Group and
              Centre for Computational Engineering,
              Technische Universität Darmstadt,
              Darmstadt, Germany \email{melina.merkel@tu-darmstadt.de}
\and R. Vázquez \at Departamento de Matem\'atica Aplicada,
              Universidade de Santiago de Compostela,
              Santiago de Compostela, Spain
              and
              Galician Centre for Mathematical Research and Technology (CITMAga),
              Santiago de Compostela, Spain \email{rafael.vazquez@usc.es}}
%
%
\maketitle


\abstract*{
In magnetostatics and eddy current problems, formulated in terms of the magnetic vector potential, the solution is not unique, 
because the addition of an irrotational function to the solution remains a valid solution. The tree-cotree decomposition is a 
gauging technique to recover uniqueness when using finite elements, which consists in considering the mesh as a graph, and building a spanning 
tree on that graph. The idea has been recently extended to isogeometric analysis, applying the construction of the spanning tree 
on the control mesh, or equivalently, on the Greville grid. In the present paper we extend the construction to hierarchical 
splines, a set of splines with multi-level structure for adaptive refinement, by constructing a spanning tree for each 
single level. Since for degree $p=1$ the spaces of finite elements and hierarchical splines coincide, the presented 
construction is also valid for quadrilateral finite element meshes with hanging nodes. To assess the correctness of the method, we 
present numerical results for Maxwell eigenvalue problem.
}
\abstract{
In magnetostatics and eddy current problems, formulated in terms of the magnetic vector potential, the solution is not unique, 
because the addition of an irrotational function to the solution remains a valid solution. The tree-cotree decomposition is a 
gauging technique to recover uniqueness when using finite elements, and it consists in considering the mesh as a graph, and building a spanning 
tree on that graph. The idea has been recently extended to isogeometric analysis, applying the construction of the spanning tree 
on the control mesh, or equivalently, on the Greville grid. In the present paper we extend the construction to hierarchical 
splines, a set of splines with multi-level structure for adaptive refinement, by constructing a spanning tree for each 
single level. Since for degree $p=1$ the spaces of finite elements and hierarchical splines coincide, the presented 
construction is also valid for quadrilateral finite element meshes with hanging nodes. To assess the correctness of the method, we 
present numerical results for the Maxwell eigenvalue problem.
}

\section{Introduction}
\label{sec:intro}
The equations for magnetostatics and magneto-quasistatics, or eddy currents,
which describe the solution of slowly varying magnetic fields \cite{Jackson_1998aa},
are commonly formulated in terms of the magnetic vector potential. For example, assuming homogeneous boundary conditions,
the magnetostatic formulation is given by 
\begin{align} \label{eq:weak}
 \int_\Omega \nu\, \curl \Afield_h  \cdot \curl \field{v}_h \operatorname{d}\Omega = \int_\Omega  \Jsrc \cdot \field{v}_h \operatorname{d}\Omega,
 \qquad \text{ for all } \field{v}_h \in V_h,
\end{align}
where $\nu$ is the magnetic reluctivity, $\Jsrc$ is an excitation current and $V_h \subset \mathbf{H}_0(\mathbf{curl};\Omega)$
is a curl-conforming discrete space. These 
formulations are usually discretized with edge finite elements, and it is known that the solution is 
unique up to an irrotational field. To recover uniqueness, one of the possible gauging
techniques, which works directly at discrete level, is tree-cotree gauging 
\cite{Albanese_1988aa,Manges_1995aa,Munteanu_2002aa}. The idea of the method is to consider the 
vertices and edges of the finite element mesh as nodes and edges of a graph, and then construct a spanning tree on this graph. 
The spanning tree is a subgraph that passes through 
every node without closing any loop, and \rv{on each connected component} it contains as many edges as the number of nodes minus one\rv{, see for example \cite[Sec.~7]{Strang_1986aa} or \cite[Def.~5.2]{Bossavit_1998aa}}.
Selecting the basis functions associated to edges in the complement of the tree, which is called 
the cotree, leads to a linear system which is uniquely solvable. The application of tree-cotree techniques 
to high order finite elements can be achieved either with hierarchical conforming elements\footnote{The construction
of hierarchichal finite elements, where new functions are added for each degree, 
is different from hierarchical B-splines, where functions are added and removed at each level.} \cite{Schoberl_2005aa}, 
applying the tree construction only to low order functions, 
or with a careful construction of the basis of edge elements \cite{Alonso-Rodriguez_2024aa}.

An alternative to high order finite elements that has gained popularity is isogeometric analysis (IGA), a 
discretization method based on splines \cite{Hughes_2005aa}, and an analogue to edge finite elements
are the curl-conforming spline spaces introduced in \cite{Buffa_2010aa}. It has been shown in \cite{Kapidani_2022aa}
that tree-cotree techniques can be easily generalized in the case of tensor-product \rv{splines by} exploiting 
the existence of commutative isomorphisms between spline spaces and finite element spaces defined in an auxiliary grid, called the Greville grid. 
Tree-cotree techniques in IGA have also been used to explicitly construct a basis of cohomology generators
for formulations in terms of the scalar magnetic potential \cite{Kapidani_2024aa}, and also in conjunction
with domain decomposition methods \cite{Mally_2025ab}.

One of the drawbacks of the spline spaces from \cite{Buffa_2010aa} is their tensor-product structure, which only
permits global refinement. To remove this constraint, spaces with local refinement capabilities have been studied
in recent years in IGA, with (truncated) hierarchical B-splines being probably the most popular ones 
\cite{Vuong_2011aa,Giannelli_2012aa}. They are based on a simple multi-level structure, where the space of each 
level is tensor-product, and the set of active functions from each level is chosen from their supports. 
As for the tensor-product case, curl-conforming spaces with hierarchical splines have been introduced in 
\cite{Evans_2020}, along with certain conditions on the hierarchical meshes to guarantee the absence of spurious 
harmonic functions \cite{Evans_2020,Shepherd_2024aa}.

While the construction of the tree to gauge the problem is easy for tensor-product splines, the extension to
hierarchical splines is not obvious. It is not possible to build a global tree across all refinement levels,
because the correspondence with finite elements in an auxiliary grid does not hold anymore. Moreover, the meshes 
of hierarchical splines are analogous to quadrilateral/hexahedral meshes with hanging nodes in finite elements,
and we are not aware of any extension of tree-cotree gauging in this setting. \rv{A different multi-level gauging technique has been presented in \cite{Hiptmair_2000aa} for tetrahedral and hexahedral finite elements. In this method, edges from the same level are gathered if they are in the interior of the same (macro)element of the previous coarser level, and simple local solves are applied on this macroelement.}
In this paper we \rv{follow a different approach, and} introduce tree-cotree techniques for spaces
of hierarchical splines in simply connected domains with full Dirichlet boundary:
exploiting the multi-level structure of hierarchical splines, a tree is built for each 
level, creating an object that we define as a multi-level tree. We present the algorithm for the construction of
the tree, and we show that it provides a valid gauge using numerical tests. It is also worth to remark that, for degree $p=1$, hierarchical spline spaces coincide with finite 
elements. Therefore, our method provides an extension of tree-cotree gauging to finite element meshes with hanging nodes.

The outline of the paper is as follows. In \Cref{sec:B-splines} we recall the main concepts to define the 
tensor-product spline spaces on a single level, and the associated Greville grid in which we build the tree. 
\Cref{sec:HB-splines} introduces the definitions of hierarchical splines, and extends the concept of Greville grids
to Greville subgrids, which will be necessary to define the multi-level tree. The construction of the tree is 
detailed in \Cref{sec:hierarchicalTC}, along with the algorithm to build it in practice. 
The correctness of the computed gauge is validated in \Cref{sec:numResults} by numerical tests, solving the Maxwell eigenvalue
problem to prove that gauging does not affect the computed eigenvalues.

\section{Splines defined on a single level}
\label{sec:B-splines}
In this section we introduce the spline de Rham sequence for dimension two, both for single-patch and 
multi-patch domains. These spaces will be used to construct the corresponding sequences of hierarchical splines
in the next section. We also present the definition of the Greville grids, that will be necessary to apply 
the tree-cotree gauging in \Cref{sec:hierarchicalTC}.

\subsection{Spline de Rham sequence in a single-patch domain}
Given two integers $n>0$ and $p>0$, and a knot vector $\gls{knotvec} = \{\xi_1, \xi_2, \ldots, \xi_{n+p+1}\}$, 
with $0 = \xi_1 \le \xi_2 \le \ldots \le \xi_{n+p+1} = 1$, the B-spline basis of degree $p$ is obtained using the 
Cox-de Boor formula, and we denote by $S_p(\Xi)$ the corresponding spanned space. 
Multivariate B-splines of degree $p$ are defined in the unit square $\gls{Omegahat} = (0,1)^2$ 
by tensor-product, from given values $n_1, n_2> 0$ and knot vectors $\Xi_1, \Xi_2$.
It is well known that a discrete de Rham sequence can be constructed using B-splines of mixed degrees $p$ and $p-1$ \cite{Buffa_2010aa,Buffa_2011aa}. 
We will denote the spaces of the sequence as \gls{Xkhat}, 
for $k = 0,1, 2$, \rv{which are defined as}
\begin{align*}
	\hat{X}^0 = &S_p(\Xi_1) \otimes S_p(\Xi_2), \\
	\hat{X}^1 = & S_{p-1}(\Xi'_1) \otimes S_p(\Xi_2) \times S_{p}(\Xi_1) \otimes S_{p-1}(\Xi'_2), \\
	\rv{\hat{X}^2 = }&\rv{S_{p-1}(\Xi'_1) \otimes S_{p-1}(\Xi'_2),}
\end{align*}
where $\Xi_i'$ is built from $\Xi_i$ removing the first and last knots. 
Moreover, we also denote by $\gls{Bkhat}$ 
the basis for $\hat{X}^k$, for $k = 0,1,2$,
which is built from standard B-splines for degree $p$ and Curry-Schoenberg splines for degree $p-1$, as in \cite{Ratnani_2012aa}. \rv{In general, when we mention degree $p$ we will refer to the degree of the first space of the sequence.} \rv{For this work, we are mainly interested in the spaces $\hat{X}^0$ and $\hat{X}^1$.}

We assume that our physical domain $\gls{Omega} \subset \mathbb{R}^2$ 
is a surface constructed through a parametrization $\gls{F}: \hat \Omega \rightarrow \Omega$ 
and, as standard in IGA, in practice it will be defined by splines or NURBS, with the same degree and continuity as $\hat X^0$. \rv{More detailed assumptions on $\gls{F}$ are given in \cite[Assumption 3.1]{Beirao-da-Veiga_2014aa}} Then, the discrete spaces \gls{Xk} in the physical domain are defined by suitable pull-backs, and in particular we have
\begin{align*}
& X^0(\Omega) = \{ u : u \circ \map \in \hat{X}^0\}, \quad \\
& X^1(\Omega) = \{ \mathbf{u} : D\mathbf{\map}^\top (\mathbf{u} \circ \map) \in \hat{X}^1\}, \\
& \rv{X^2(\Omega) = \{ u : \det(D\mathbf{\map}) (u \circ \map) \in \hat{X}^2\},}
\end{align*}
where $\gls{DF}$ 
is the Jacobian of the parametrization. The bases $\gls{Bk}$ 
are defined applying the same pull-backs to the basis functions in the parametric domain.

\subsection{Spline de Rham sequence in multi-patch domains}
More complex domains can be defined using multi-patch techniques. We assume that \rv{$\Omega$ is a connected Lipschitz domain} formed by several patches, 
$\overline{\Omega} = \bigcup_{i=1}^{N_p} \overline{\gls{Omegai}}$, 
each patch parametrized in the form $\gls{Fi}: \hat \Omega \rightarrow \Omega^{(i)}$ 
and satisfying that $\Omega^{(i)} \cap \Omega^{(j)} = \emptyset$ for $i \not = j$. Moreover, we assume that $\partial \Omega^{(i)} \cap \partial \Omega^{(j)}$ 
corresponds to either a vertex or a full edge for each patch. Assuming that the two parametrizations at the interfaces 
coincide up to an affine transformation, and that the same holds for their corresponding knot vectors, 
one can define discrete conforming spaces in the multi-patch domain as
\begin{align*}
& X^0(\Omega) = \{ u \in H^1(\Omega) : u|_{\Omega^{(i)}} \in X^0(\Omega^{(i)})\}, \, \\
& X^1(\Omega) = \{ \mathbf{u} \in \Hcurl: \mathbf{u}|_{\Omega^{(i)}} \in X^1(\Omega^{(i)})\}, \\
& \rv{X^2(\Omega) = \{ u \in L^2(\Omega) : u|_{\Omega^{(i)}} \in X^2(\Omega^{(i)})\} }.
\end{align*}
In this situation, there is a one-to-one correspondence between functions of neighboring patches that do not vanish on the interface, and they can be identified up to orientation, as it is done in finite elements. We refer the reader to \cite[Sect.~4.4]{Buffa_2013aa} for the details on how to define the spaces in the multi-patch setting.
By abuse of notation, we will denote in the same way the bases of the single-patch and multi-patch cases, namely $\gls{Bk}$.

\subsection{The Greville grid and auxiliary finite element spaces} \label{sec:Greville}
In the univariate setting, each basis function is associated to a Greville abscissa defined by the knot average $\gls{Pi} = (\xi_{i+1} + \ldots + \xi_{i+p}) / {p}$, 
for $i = 1, \ldots, n$. In the unit square, we can associate to each basis function in $\hat \basis^0$ a Greville point, 
defined as the Cartesian product of univariate points, and map them through $\map$ to obtain the corresponding Greville 
points for $\basis^0$, the basis functions in the physical domain. By connecting the Greville points we obtain a Cartesian 
grid in the parametric domain $\hat \Omega$, and a structured grid of bilinear quadrilaterals that approximate the physical 
domain $\Omega$. We will refer to this as the Greville grid, and it will be denoted by $G$. By construction, the basis 
functions in $\basis^0$ are associated to the nodes of the Greville grid, the basis functions in $\basis^1$ are 
associated to its edges, \rv{and the basis functions in $\basis^2$ are associated to its cells}. Moreover, we can define in the grid $\gls{G}$ 
the discrete spaces of lowest order nodal elements and lowest order finite edge elements, that we denote by $\glslink{Zk}{Z^0}$ and $\glslink{Zk}{Z^1}$. 
Then, identifying each basis function of the spline spaces with the corresponding basis functions of these auxiliary finite element spaces, there exist linear isomorphisms $\glslink{Ik}{I^0}, \glslink{Ik}{I^1}$ 
such that the following diagram commutes \cite{Buffa_2013aa}
\begin{equation} \label{eq:comm}
	\begin{CD}
		X^0(\Omega) @>\grad>> X^1(\Omega) \\
		@V{I^0}VV  @V{I^1}VV \\
		Z^0 @>\grad>> Z^1.
	\end{CD}
\end{equation}

\begin{figure}[th]
  \includegraphics[width=0.45\textwidth]{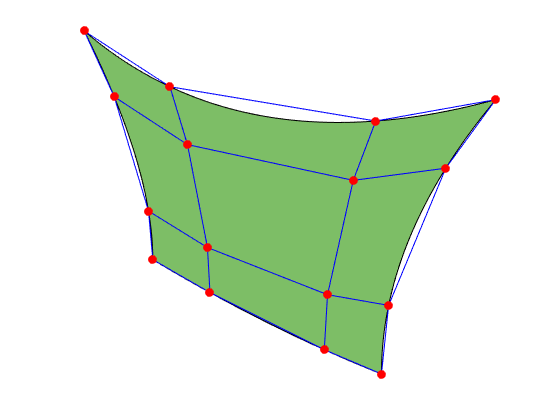}
  \includegraphics[width=0.45\textwidth]{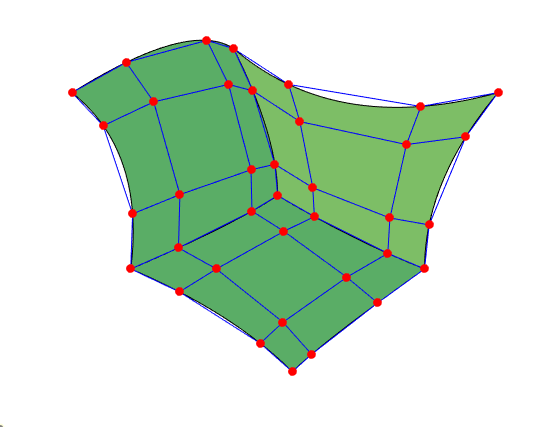}
  \caption{Example of two-dimensional Greville grids for the single-patch (left) and multi-patch case (right), \rv{for degree $p=3$ and knot vector $\gls{knotvec} = \{0, 0, 0, 0, 1, 1, 1, 1\}$ in all patches}.}
  \label{fig:greville_grids}
\end{figure}

In the multi-patch setting, due to the assumptions on how the patches are defined, the Greville points at the
interface between two neighboring patches coincide. Therefore, it is possible to associate a Greville grid to the 
spline spaces also in this more general setting, with the difference that the grid is structured on each patch, but 
globally unstructured. The finite element spaces and the isomorphisms can also be defined in the multi-patch setting, and
abusing notation, we will denote them as in the single-patch case. Examples of Greville grids are shown in \Cref{fig:greville_grids}.
To simplify the presentation, in the following the examples will be given for single-patch domains, and the visualization will be done 
for the Greville grid of the parametric domain.

\begin{remark}
	The Greville grid and the lowest order finite element spaces are auxiliary tools for the tree-cotree algorithm 
	that will be introduced afterwards. Since the tree-cotree algorithm is based on topological properties, \rv{it will work for any other 
	mesh with the same number of nodes and edges, and the same connectivity.}
\end{remark}

\section{Hierarchical B-splines}
\label{sec:HB-splines}
Hierarchical splines are defined using a multi-level structure. First, we define standard spaces of splines on each level,
as in the previous section, that can be either single-patch or multi-patch. Then, by choosing the regions where we want to 
apply local refinement, we apply a selection algorithm in which basis functions are activated according to the relation of their supports with these regions
\cite{Kraft_1997aa,Vuong_2011aa}.
While the focus will be on the first two spaces of the de Rham sequence, for completeness we introduce the notation for the whole sequence,
introduced in \cite{Evans_2020}.

From now on we will assume that the spaces have vanishing boundary conditions, i.e., we will work with 
$X^0(\Omega) \cap H^1_0(\Omega)$ and $X^1(\Omega) \cap \mathbf{H}_0(\mathbf{curl};\Omega)$. To simplify notation, we will
use the same notation $X^0(\Omega)$ and $X^1(\Omega)$ for the spaces with vanishing boundary conditions. \rv{We also recall that we are assuming that the domain $\Omega$ has trivial topology.}

\subsection{Definition of hierarchical splines} We start defining for $k = 0,1,2$ a sequence of nested spline spaces
\[
X^k_0 \subset X^k_1 \subset \ldots \subset X^k_{L-1} \subset X^k_L,
\]
where for every level $\ell$, the spaces $\gls{Xkl}$ 
for $k>0$ are defined from $X^0_\ell$, and they form a de Rham sequence  as in the previous section.
We denote the basis for each of these spaces by $\gls{Bkl}$, 
for each level $\gls{l} = 0, \ldots, L$. 
Moreover, we notice that at each level we have a Greville grid $\gls{Gl}$. 

To define hierarchical B-splines, we also need a nested sequence of closed subdomains
\begin{equation*} 
\overline{\Omega} = \glslink{Omegal}{\Omega_0} \supseteq \glslink{Omegal}{\Omega_1} \supseteq \ldots \supseteq \glslink{Omegal}{\Omega_L} \supseteq \glslink{Omegal}{\Omega_{L+1}} = \emptyset ,
\end{equation*}
where each subdomain is defined as the union of supports of functions from the previous level, namely
\begin{equation} \label{eq:union_of_supports}
	\Omega_{\ell+1} = \bigcup_{\beta \in \mathcal{S}_{\ell}} \overline{\supp} \beta, \quad \text{for some } \mathcal{S}_{\ell} \subset \basis^2_{\ell}, \, \ell = 0, \ldots, L-1.
\end{equation}

Let us define, for $\ell = 0, \ldots, L$ and $\ell' \in \left\{\ell, \ell+1\right\}$, 
the subset of basis functions with (open) support completely contained in $\Omega_{\ell'}$, namely
\[
\gls{Bkll} = \{ \beta \in \basis^k_\ell : \supp \beta \subset \Omega_{\ell'}\}.
\]
Then, the basis of hierarchical B-splines $\gls{HkL}$ 
is defined, for $k = 0, 1, 2$, by applying the following recursive algorithm:
\begin{enumerate}
	\item $\hbasis^k_0 = \basis^k_0$,
	\item $\hbasis^k_{\ell+1} = (\gls{Hkl} \setminus \basis^k_{\ell,\ell+1}) \cup \basis^k_{\ell+1, \ell+1}$, for $\ell = 0, \ldots, L-1$.
\end{enumerate}
That is, we start from the basis of tensor-product B-splines for the coarsest space, and at each step we remove coarse basis functions
with support contained in $\Omega_{\ell+1}$, and add fine basis functions with support contained in the same subdomain. It has been 
proved in \cite{Vuong_2011aa} (see also \cite{Evans_2020}) that the functions in $\hbasis^k_{\ell}$ are linearly independent for any $\ell$. 
We will denote by $\gls{Wkl} = \mathrm{span} \{\hbasis^k_{\ell}\}$ the space of hierarchical B-splines up to level $\ell$.

We also define the set of active basis functions of level $\ell$, and the set of deactivated functions of level $\ell$, respectively as 
\begin{equation} \label{eq:active_functions}
\gls{Akl} = (\basis^k_{\ell,\ell} \setminus \basis^k_{\ell,\ell+1}),
\qquad
\gls{Dkl} = \basis^k_{\ell,\ell+1}.
\end{equation}
It is easy to see, due to the nestedness of the subdomains, that 
$\hbasis^k_{\ell} = (\bigcup_{\ell' = 0}^{\ell-1} \act^k_{\ell'}) \cup \basis^k_{\ell,\ell} = (\bigcup_{\ell' = 0}^{\ell} \act^k_{\ell'}) \cup \deactivated^k_{\ell}$, and as a consequence
$\hbasis^k_L = \bigcup_{\ell = 0}^{L} \act^k_\ell$, because $\Omega_{L+1} = \emptyset$. 

\subsection{The Greville subgrids and associated finite element spaces}
We now introduce the Greville subgrids, that will play a prominent role in the tree-cotree decomposition.
Following \cite{Evans_2020}, we denote by $\gls{Gll}$, 
with $\ell' \in \left\{\ell, \ell+1\right\}$, the Greville subgrid formed by the cells \rv{of the Greville grid $G_\ell$} associated to basis functions 
in $\basis^2_{\ell,\ell'}$, that is, functions of level $\ell$ with support completely contained in $\Omega_{\ell'}$. \rv{By abuse of notation, we will also denote by $G_{\ell,\ell'}$ the domain covered by (the interior of the closure of) the union of the cells in the 
Greville subgrid.}
We note that, \rv{since we have assumed vanishing boundary conditions}, a basis function belongs to $\basis^k_{\ell,\ell'}$ \rv{for $k = 0,1,2$} if and only if it is associated to an internal entity of the subgrid $G_{\ell,\ell'}$. \rv{Indeed, from the boundary conditions, functions associated to boundary entities of $G_\ell$ do not even belong to the basis $\basis^k_{\ell} \supset \basis^k_{\ell,\ell'}$. Moreover, for $k=0,1$, if a function in $\basis^k_{\ell,\ell'}$ is associated to an internal entity of $G_\ell$, then the functions associated to adjacent cells must belong to $\basis^2_{\ell,\ell'}$, and therefore it is also an internal entity of $G_{\ell,\ell'}$. The assumption on homogeneous Dirichlet boundary conditions allows to treat in the same manner all the boundary entities, simplifying the presentation of the method. While mixed boundary conditions can also be considered, this would increase the complexity of the notation and the presentation because functions associated to Neumann boundaries should be treated as internal.}

\begin{figure}
\begin{subfigure}{0.32\textwidth}
\includegraphics[width=\textwidth,trim=14cm 1cm 12cm 1cm, clip]{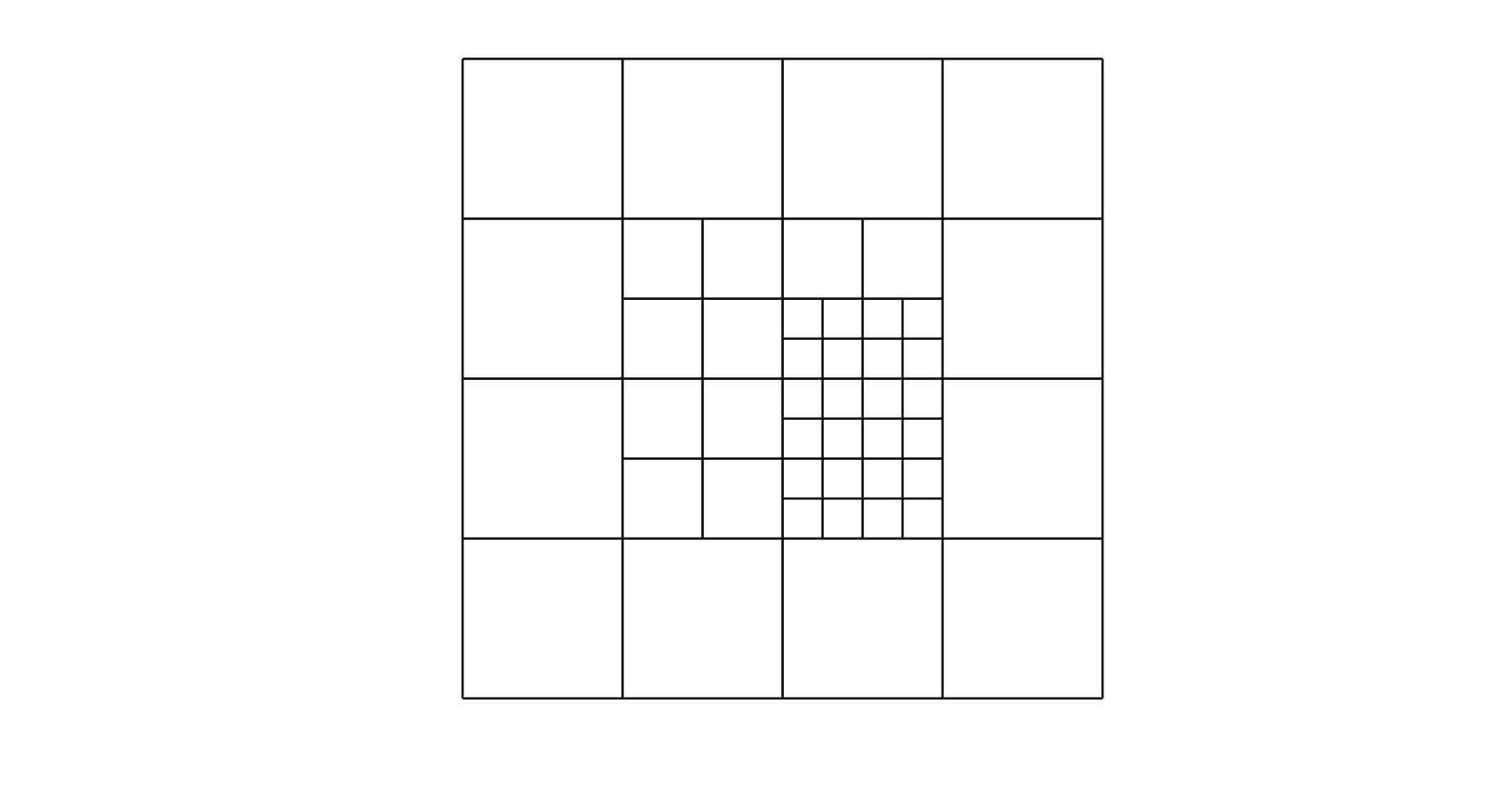}
\caption{Hierarchical mesh, $p=1$}
\end{subfigure}
\begin{subfigure}{0.32\textwidth}
\includegraphics[width=\textwidth,trim=14cm 1cm 12cm 1cm, clip]{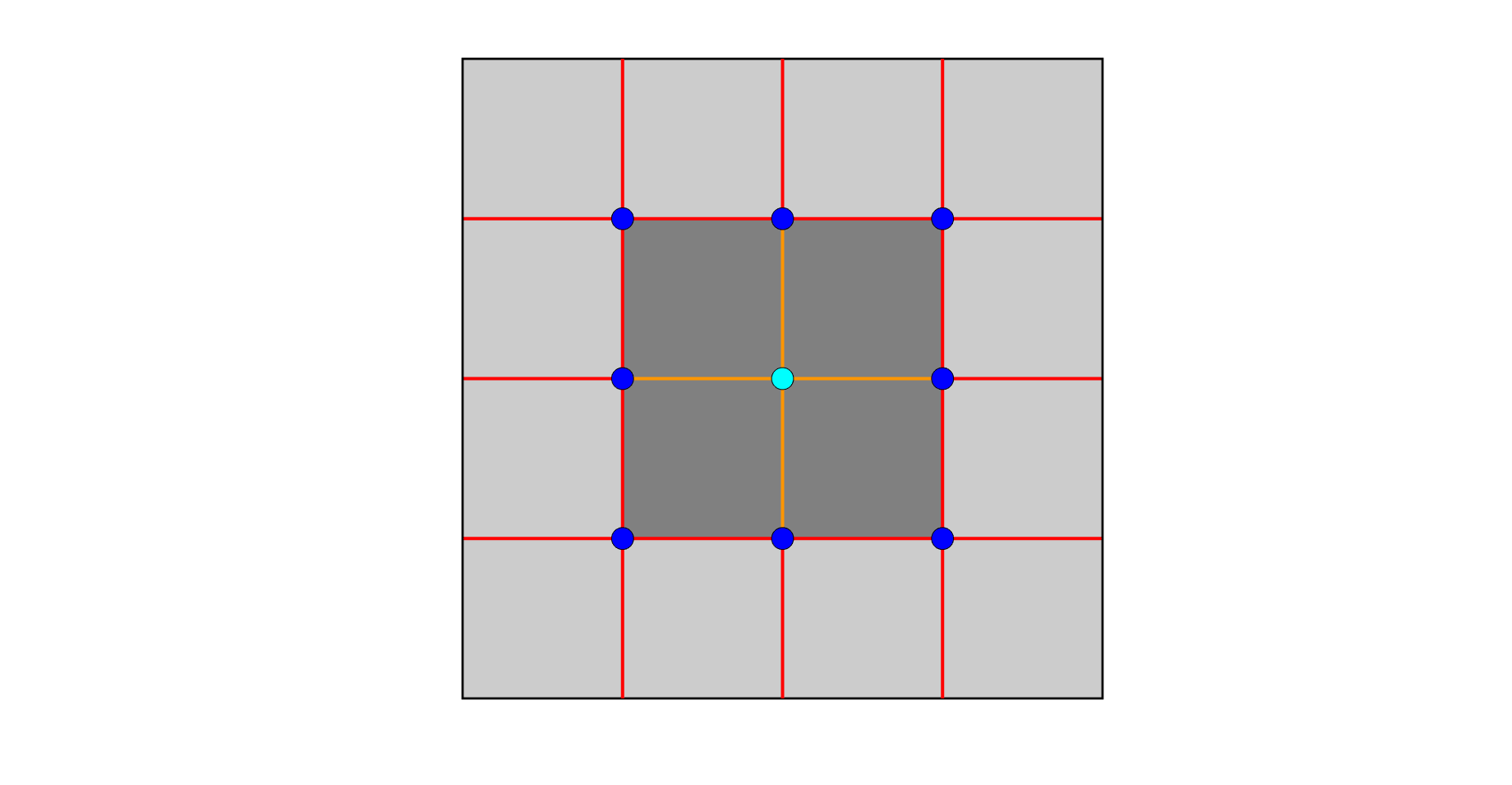}
\caption{Greville subgrid of level 0}
\end{subfigure}
\begin{subfigure}{0.32\textwidth}
\includegraphics[width=\textwidth,trim=14cm 1cm 12cm 1cm, clip]{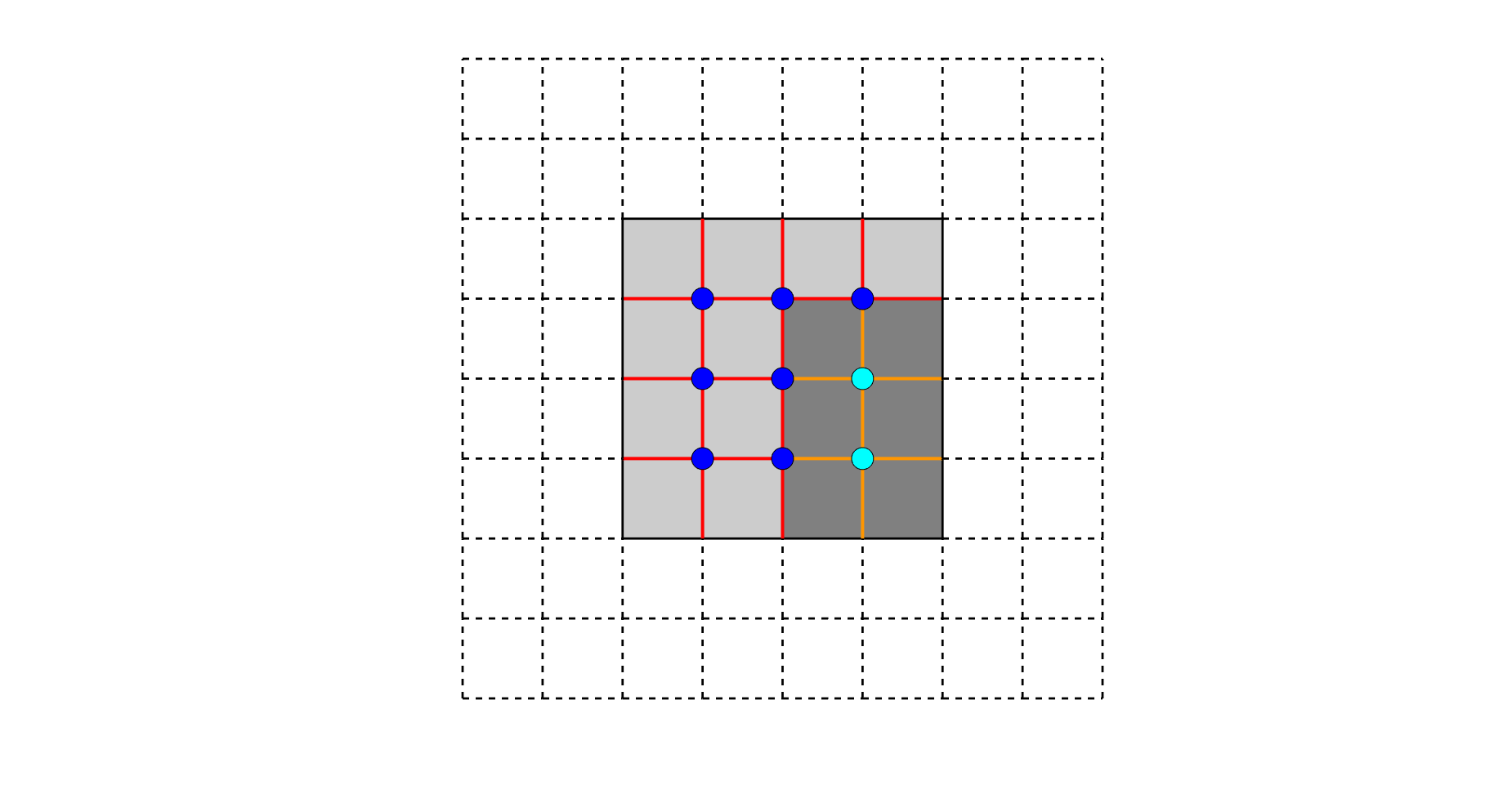}	
\caption{Greville subgrid of level 1}
\end{subfigure}

\begin{subfigure}{0.32\textwidth}
\includegraphics[width=\textwidth,trim=14cm 1cm 12cm 1cm, clip]{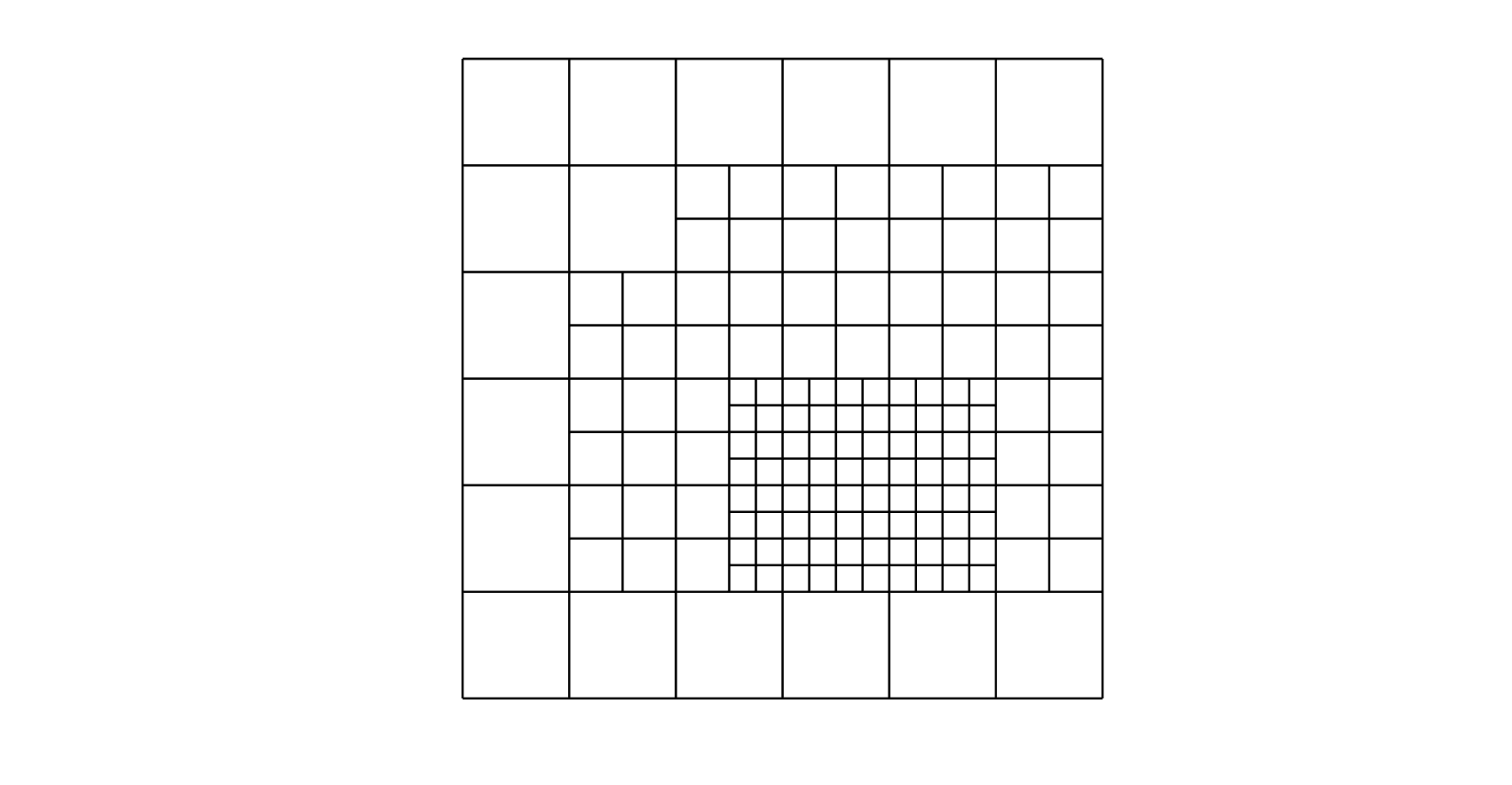}
\caption{Hierarchical mesh, $p=2$}
\end{subfigure}
\begin{subfigure}{0.32\textwidth}
\includegraphics[width=\textwidth,trim=14cm 1cm 12cm 1cm, clip]{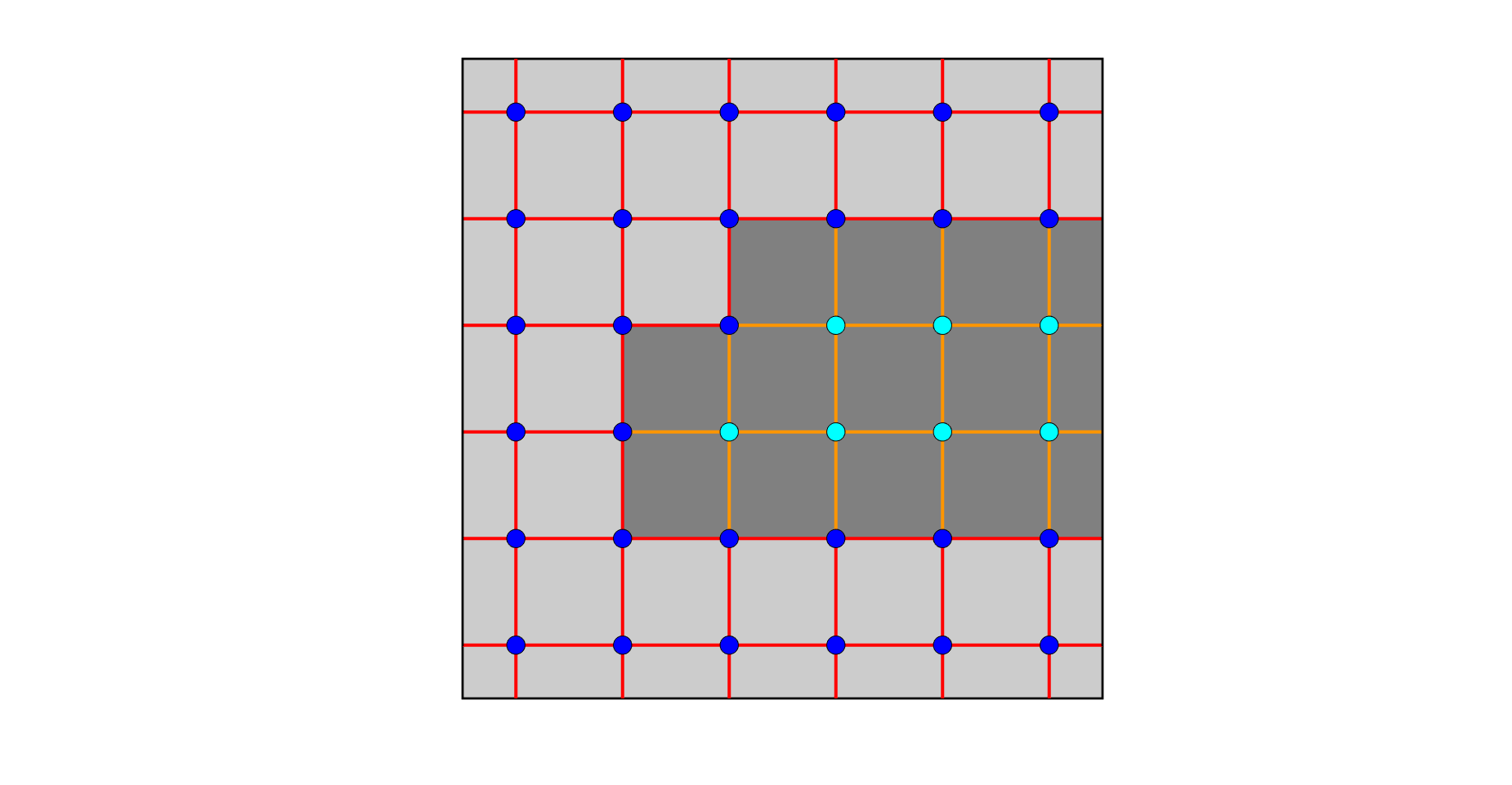}
\caption{Greville subgrid of level 0}
\end{subfigure}
\begin{subfigure}{0.32\textwidth}
\includegraphics[width=\textwidth,trim=14cm 1cm 12cm 1cm, clip]{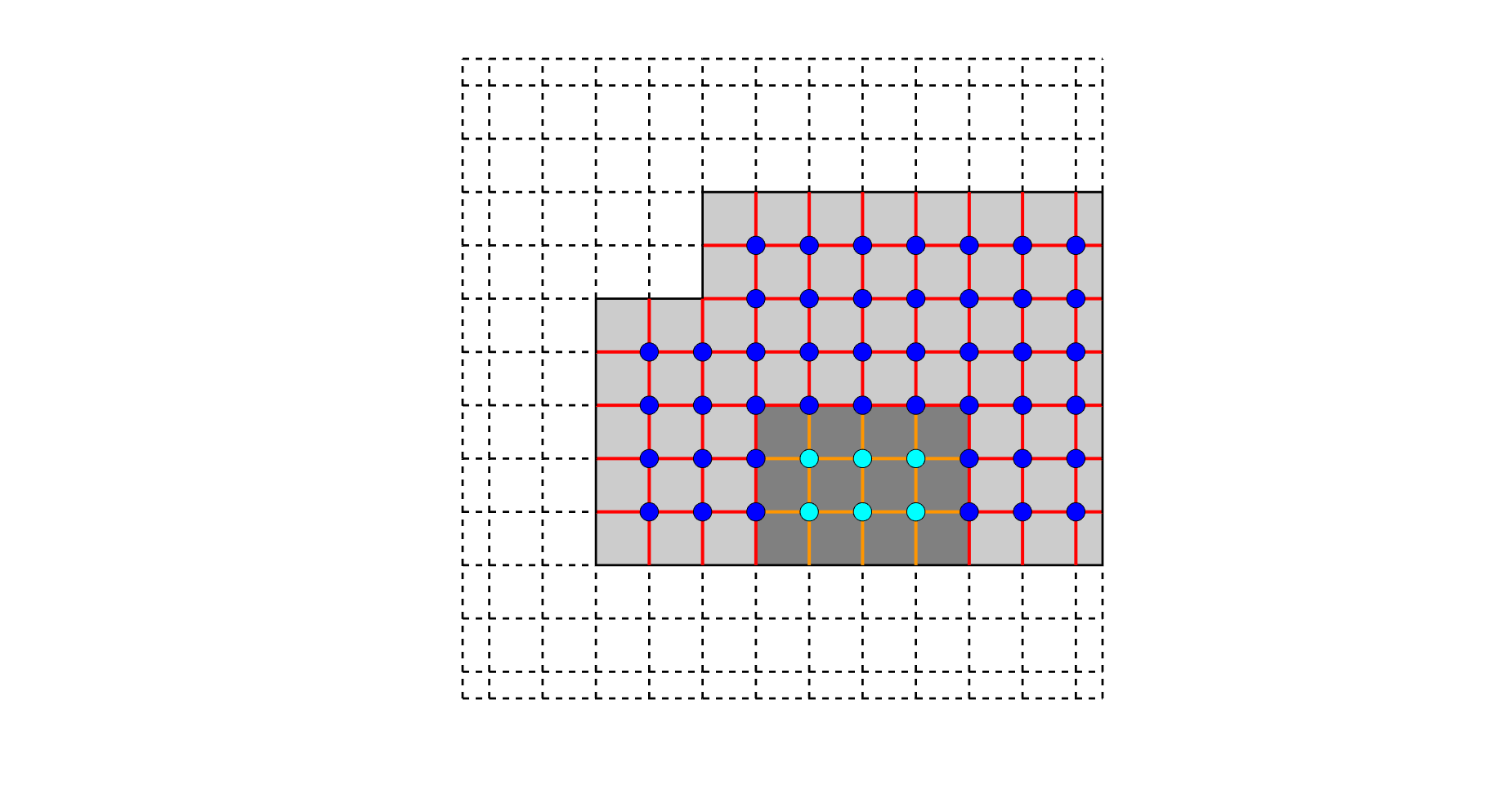}	
\caption{Greville subgrid of level 1}
\end{subfigure}
\caption{Example for degrees $p=1$ (top row) and $p=2$ (bottom row) of hierarchical meshes with three levels (left), 
and the corresponding Greville subgrids for levels $\ell=0$ (center) and $\ell=1$ (right). The Greville subgrid
$G_{\ell,\ell}$ corresponds to all the highlighted gray region, and the Greville subgrid $G_{\ell,\ell+1}$ corresponds 
to the dark gray region. The entities of the Greville subgrids associated to active basis functions are marked with blue
vertices ($\act^0_\ell$) and red edges ($\act^1_\ell$), while those associated to deactived functions are marked with 
cyan vertices ($\deactivated^0_\ell$) and orange edges ($\deactivated^1_\ell$).
 Note that all the internal entities of $G_{\ell,\ell+1}$ correspond to deactivated functions.}
\label{fig:greville_subgrids}
\end{figure}

In \Cref{fig:greville_subgrids} we present two examples of Greville subgrids for hierarchical splines of degrees $p=1$ and $p=2$.
First of all, note that for degree $p=1$ there is a clear correspondence between the hierarchical mesh and the Greville 
subgrids, but this is lost when increasing the degree, due to the definition of the Greville points as knot averages.
In fact, due to the difference of the knot averages between levels, the domains $G_{\ell,\ell+1}$ and 
$G_{\ell+1,\ell+1}$ are in general not equal, even if both correspond to functions with support contained in 
$\Omega_{\ell+1}$. These geometric differences between the grids will not be important for the tree-cotree decomposition, 
because the construction of the tree will be based on the topology of the mesh, \rv{and the Greville subgrids will have the same topology under the assumptions of \Cref{sec:assumptions}}. 
We also remark that the active \rv{Greville} entities associated to active functions in $\act^k_{\ell}$ of the hierarchical basis are contained 
in the closure of the subdomain $G_{\ell,\ell} \setminus G_{\ell,\ell+1}$, while those associated to deactivated functions
are contained in the interior of $G_{\ell,\ell+1}$.

\rv{
\begin{remark}
	As mentioned above, for degree $p=1$ there is a correspondence between the cells of the Greville subgrids and the elements of the hierarchical mesh. This can be exploited to prove that, as in the tensor-product case, for $p=1$ the spaces of hierarchical splines coincide with finite elements defined on the hierarchical mesh (with hanging nodes), and for all the spaces of the de Rham sequence. The result is trivial for $k=2$, since the functions are piecewise constants. For $k=0,1$, it can be proved by the fact that the number of basis functions is the same, they are piecewise polynomials of the same degree and continuity, and in this case hierarchical splines span the whole space of piecewise polynomials \cite{Mokris2014}. While the result is probably true also for $p>1$ with low continuity, the proof in this case is not trivial, because we lack the correspondence of the meshes.
\end{remark}
}

\subsection{Constraints for exactness of the sequence} \label{sec:assumptions}

\rv{In the single level construction, and under the assumption that the domain has trivial topology, the spline spaces form an exact sequence. Instead, although}
the construction of the de Rham sequence of hierarchical splines is possible for any hierarchical mesh,
there are situations in which the spaces of hierarchical splines $W^k(\Omega)$ do not form an exact sequence, 
and this may cause the appearance of spurious results. A necessary and sufficient condition for exactness\rv{, which implies that $G_{\ell,\ell} $ and $G_{\ell,\ell+1}$ have the same topology}, was first 
identified in \cite[Thm.~5.5]{Evans_2020}\rv{, but} in the same paper it was shown that there exist situations in which
the sequence is exact and spurious results still appear. For this reason the authors introduced a stronger sufficient condition 
for exactness based on the structure of the hierarchical mesh, \rv{which removed spurious results in all the numerical tests, and} which reads as follows
\begin{assumption} \label{assum:Evans_et_al}
For each level $\ell = 0, \ldots, L$ and for any basis function $\beta \in \mathcal{B}^k_\ell$, with $k = 0,1,2$, the region
$\supp{\beta} \cap (\Omega \setminus \overline{\Omega_{\ell+1}})$ is connected and simply connected.
\end{assumption}
Another sufficient condition, which covers different cases, was presented in \cite{Shepherd_2024aa}, with the main advantage
that it can be generalized to arbitrary dimension. We also remark that a refinement algorithm that automatically satisfies
this condition was recently introduced, in the two-dimensional case, in \cite{Cabanas_2025aa}.
\begin{assumption}\label{assum:exactness}
	Let us assume that in \eqref{eq:union_of_supports} we have $\mathcal{S}_\ell \subset \basis^0_\ell$.
	Moreover, let $\beta_1, \beta_2 \in \basis^0_{\ell,\ell+1}$. If they share an $(\ell+1)$-intersection, as in 
	\cite[Def.~4.1]{Shepherd_2024aa} with $n=2$, they must be connected through a shortest chain of 
	functions in $\basis^0_{\ell,\ell+1}$, as in \cite[Def.~4.3]{Shepherd_2024aa}.
\end{assumption}
The condition consists of two parts: first, we assume that the refined region is the union of support of functions in
$\basis^0_\ell$ instead of $\basis^2_\ell$, that is, they have degree $p$ instead of degree $p-1$. Second, whenever two 
refined functions are too close, there must be other refined functions in between that connect them. We remark that the condition 
appearing in \cite{Shepherd_2024aa}, to which we refer for the details, is a bit more complex, because it covers situations
that do not occur in the two-dimensional case.

The analysis of these conditions is beyond the focus of this paper. Therefore, without further details, from now on we
will assume that the hierarchical meshes satisfy either \Cref{assum:Evans_et_al} or \Cref{assum:exactness}.

\section{Tree-cotree decomposition for hierarchical B-splines}
\label{sec:hierarchicalTC}

As already mentioned in the introduction, the solution of the magnetostatic problem \eqref{eq:weak} is not unique,
and a similar issue arises in the formulation of the magneto-quasistatic problem. One of the techniques to recover
uniqueness for lowest order finite elements is tree-cotree gauging. 
In this section we will start presenting the basics of the tree-cotree decomposition for finite elements, and its extension 
to splines of one single level, introduced for the first time in \cite{Kapidani_2022aa}. Then, we will present the extension to 
multilevel hierarchical splines, including the algorithm to construct what we call a multi-level tree.

\subsection{Tree-cotree gauging for splines of a single level}
In order to apply tree-cotree gauging for low order finite elements, as explained in 
\cite{Albanese_1988aa,Bossavit_1998aa,Manges_1995aa,Munteanu_2002aa}
the idea is to consider the finite element mesh as a graph $\gls{graph} = (\gls{N},\gls{E})$, where the nodes $\nodeSet$ and edges 
$\edgeSet$ of the graph correspond to the vertices and edges of the mesh. Then, using an algorithm from graph theory, one builds a
spanning tree $\gls{tree}$ on the graph. A spanning tree is a subgraph that passes through every node and that does 
not contain any cycle, and therefore, for a connected graph, it contains as many edges as the number of nodes minus one. 
Tree-cotree gauging is then applied by removing from the linear system the basis functions 
associated to the edges of the tree, or in other words, by keeping only those associated to the complement of the 
tree, which is called the cotree. It is important to remark that, in order to take into account Dirichlet boundary 
conditions, the tree must be grown starting from the boundary. This kind of behavior can be obtained by constructing 
a minimum spanning tree, for instance using Kruskal's algorithm \cite{Kruskal_1956aa}, and applying a lower weight 
to those edges associated to functions on which we impose Dirichlet boundary conditions. 
This procedure is also equivalent to collapsing each connected component of the Dirichlet boundary into a single node 
\cite{Rapetti_2022}.

The application of tree-cotree gauging for splines of a single level follows the same idea \cite{Kapidani_2022aa}. 
As it was explained in \Cref{sec:Greville}, the basis functions of the curl-conforming space are associated to 
the edges of the Greville grid. It is then sufficient to consider the Greville grid as a graph, and build the spanning 
tree on the Greville grid, applying exactly the same algorithms as in finite elements. Thanks to the isomorphisms between 
the spline spaces and the auxiliary finite element spaces of the Greville grid \eqref{eq:comm}, this will provide a valid gauge for the 
splines, by simply removing from the linear system the indices corresponding to basis functions associated to edges of the tree.
We note that the construction is valid both in the single-patch or in the multi-patch case. 

\subsection{Tree-cotree gauging for hierarchical splines} \label{sec:HTC_trivial}
The construction of the tree for splines on a single level is very simple, due to the existence of commutative 
isomorphisms with the auxiliary finite element spaces in the Greville grid. On the contrary, the extension to 
hierarchical splines is not obvious. In fact, it is not known if such isomorphisms exist for hierarchical splines, 
and even if that was the case, the corresponding
Greville mesh across all levels would contain hanging nodes, 
as in the meshes of \Cref{fig:greville_subgrids}. 
The presence of hanging nodes,
that have no functions attached to them, prevents an easy interpretation of the mesh as a graph on which one 
can build the spanning tree, and in fact we are not aware of any application of tree-cotree gauging on finite element
meshes with hanging nodes. 
Therefore, a new tree construction for the correct gauging of hierarchically refined B-splines is necessary and will be outlined in the following.

It is clear from \eqref{eq:active_functions} that the basis functions appearing in the hierarchical construction are 
related to the Greville subgrids. Indeed, the functions in $\basis^k_{\ell,\ell}$ and $\basis^k_{\ell,\ell+1}$
respectively correspond to the interior entities of $G_{\ell,\ell}$ and $G_{\ell,\ell+1}$. In order to obtain a valid gauge 
we will build a spanning tree on the Greville subgrid of each level, $G_{\ell,\ell}$, with some constraints to ensure that 
the space of gradient fields of functions in $W^{0}_{L}$ is correctly removed from the solution space $W^{1}_{L}$.

We start considering for each level $\ell$ the Greville subgrid $G_{\ell,\ell}$, and define the graph 
$\gls{graphl}=(\gls{Nl},\gls{El})$, 
consisting of nodes and edges $\nodeSet_{\ell}$ and $\edgeSet_{\ell}$, respectively.
Note that the graph depends on the choice of the refined region $\Omega_\ell$, and it is not necessarily connected.
%
Both the set of nodes and edges can be divided into entities in the interior of $G_{\ell,\ell}$ and on its boundary $\partial G_{\ell,\ell}$, 
namely
\begin{align*}
 \gls{Nl} &= \gls{Nintl} \cup \gls{Ndl}, & \gls{El} &= \gls{Eintl} \cup \gls{Edl}. 
\end{align*}
Due to the definition of active functions, and since we are assuming vanishing boundary conditions, 
the boundary entities always correspond to non-active functions. 
Regarding the interior entities, they can be further divided into those associated to active basis functions in 
$\act^{k}_{\ell}$ and deactivated basis functions in $\deactivated^{k}_{\ell}$, with $k=0$ for the nodes and $k=1$ for the edges, namely
\begin{align*}
 \gls{Nintl} &= \gls{NAl} \cup \gls{NDl}, &
 \gls{Eintl} &= \gls{EAl} \cup \gls{EDl}. 
\end{align*}
\rv{Note that the nodes and edges associated to deactivated functions, \gls{NDl} and \gls{EDl}, respectively correspond to the interior nodes and edges of $G_{\ell,\ell+1}$.}

\begin{figure}
 \tikzsetnextfilename{square_8x8_deg1_level0}
  \begin{subfigure}{0.24\textwidth}
  \resizebox{\linewidth}{!}{
  \includegraphics{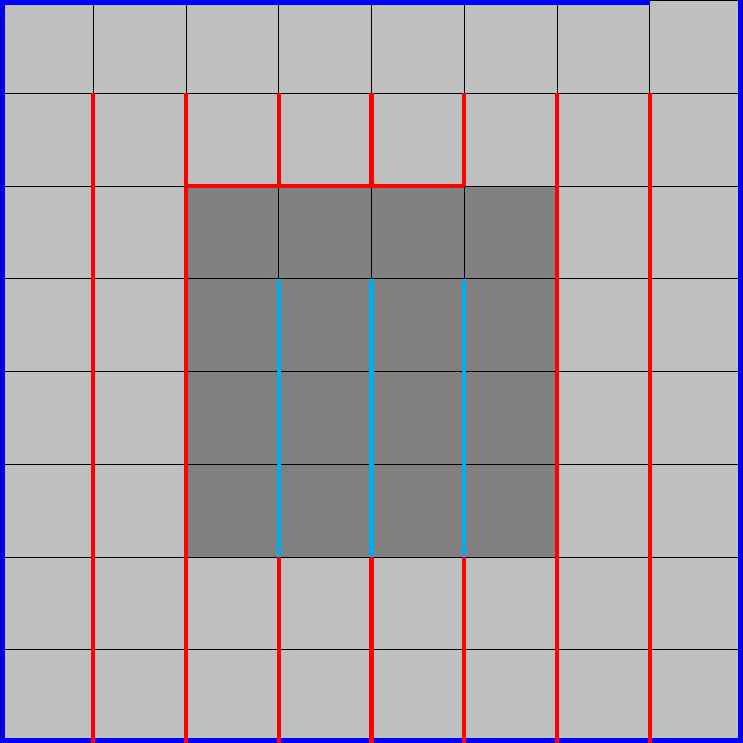}}
  \caption{Tree of level 0}
 \end{subfigure}
 \tikzsetnextfilename{square_8x8_deg1_level1}
 \begin{subfigure}{0.24\textwidth}
  \resizebox{\linewidth}{!}{
  \includegraphics{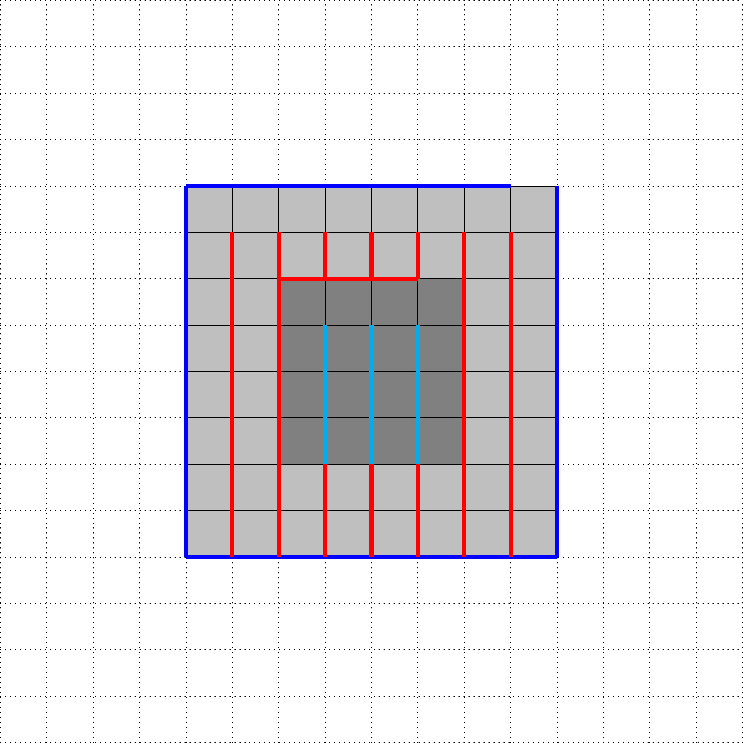}}
  \caption{Tree of level 1}
 \end{subfigure}
 \tikzsetnextfilename{square_8x8_deg1_level2}
 \begin{subfigure}{0.24\textwidth}
  \resizebox{\linewidth}{!}{
  \includegraphics{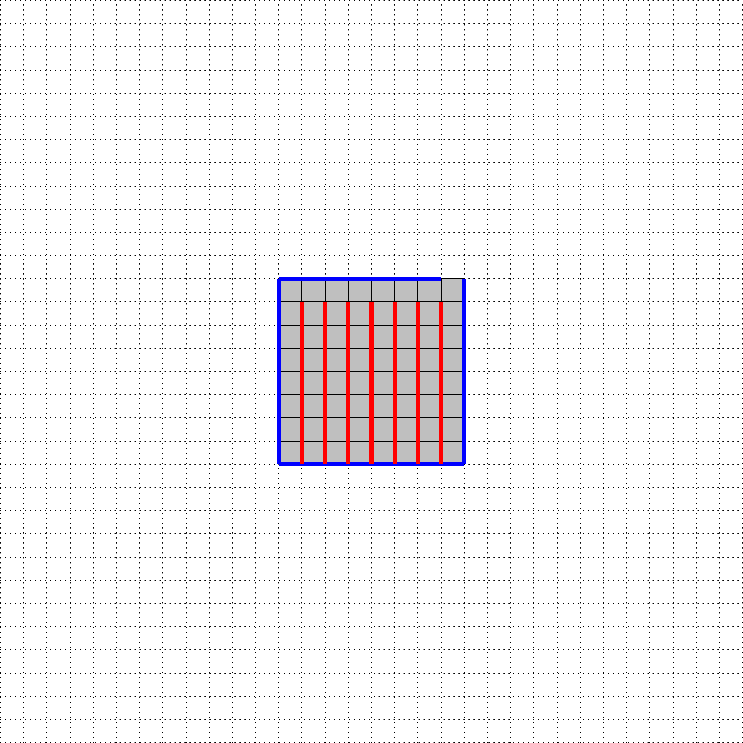}}
  \caption{Tree of level 2}
 \end{subfigure}
 \tikzsetnextfilename{square_8x8_deg1_allLevels}
 \begin{subfigure}{0.24\textwidth}
  \resizebox{\linewidth}{!}{
  \includegraphics{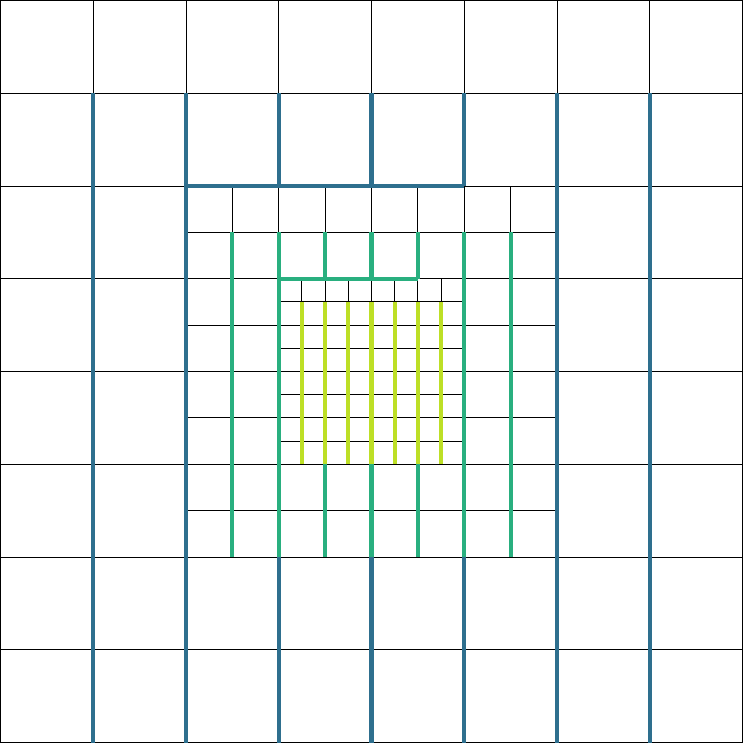}}
  \caption{Multi-level tree}
 \end{subfigure}
 \caption{Illustration of the trees built on each level for degree $p=1$, for a three-level mesh. The first three figures show 
 the tree of each single level, where the whole gray region corresponds to $G_{\ell,\ell}$, and the dark gray region corresponds
 to $G_{\ell,\ell+1}$. At each level, we first build the tree on the boundary of $G_{\ell,\ell}$ (blue edges), then on the active
 edges (red edges), and finally on the deactivated edges (cyan edges). The right figure (d) shows the multi-level tree, combining
 active edges of the three levels.}
 \label{fig:hier_TC}
\end{figure}

To achieve a correct gauging condition we propose the following approach.
On each level $\ell$ we construct a tree $\gls{Tl}$ 
using the following three steps:
\begin{enumerate}
\item Construct a spanning tree on the boundary subgraph $\gls{graphdl} = (\nodesboundary,\edgesboundary)$. 
\item Grow the tree to the interior of $G_{\ell,\ell}$ using active edges from $\edgesactive$.
\item Grow the tree to the remaining nodes using deactivated edges from $\edgesdeactivated$.
\end{enumerate}
That is, we first build a spanning tree on the boundary of $G_{\ell,\ell}$ (or more precisely, on each connected component of the boundary of $G_{\ell,\ell}$) using the boundary edges of the Greville subgrid, 
independently of whether they are associated to the Dirichlet boundary or not. Then, we grow the tree on active edges,
which \rv{are the remaining edges} in the closure of $G_{\ell,\ell} \setminus G_{\ell,\ell+1}$. Finally, 
the tree is grown on the region of deactivated functions, \rv{the interior of} $G_{\ell,\ell+1}$, only if there are nodes that have not been reached\rv{, although this last step could be removed (see \Cref{rem:deactivated_not_needed})}. 
An illustration of the procedure is shown in \Cref{fig:hier_TC} for a simple mesh and degree $p=1$.

With some abuse of notation, we will denote by $\gls{TAl}$ and $\gls{TDl}$ the restriction of the tree $\tree_\ell$ 
to edges associated to active and deactivated functions, respectively. 
The union of the active tree parts then defines the ``multi-level tree''
\begin{align*}
 \gls{ML} = \bigcup_{\ell=0}^{L} \treeactive.
\end{align*}
Since the multi-level tree consists of active edges, associated to active functions, we can define the multi-level cotree
as the set of complementary active edges (or functions), and employ them for the decomposition of the linear system 
as in the single level case. We also show in \Cref{fig:hier_TC} a visualization of the multi-level tree. Note that the plot of the 
multi-level tree is not clear for degree greater than one, because, in the Greville submeshes, active edges of different levels 
would overlap.

In practice, the construction of the multi-level tree can be achieved working on each level independently. To build the
tree on each level respecting the three steps, one can employ an algorithm to find a minimum spanning tree on weighted graphs,
such as Kruskal's algorithm \cite{Kruskal_1956aa}. The construction of the multi-level tree with such an algorithm is 
summarized in \Cref{alg:buildtree}. \mm{Note that, while we use the weights $w^\partial =1, w^\act =2, w^\deactivated=3$ any choice of weights which fulfills the relation $w^\partial < w^\act < w^\deactivated$ is possible, \rv{also with non-constant weights for each subset of edges.}}

\begin{algorithm}
\caption{Construct multi-level tree $\multileveltree_{L}$}
\label{alg:buildtree}
\begin{algorithmic}
\STATE{$\multileveltree_{L} \leftarrow \{\}$}
\FOR{$\ell=0$ to $L$}
\STATE{Set weights of edges in $\edgesboundary$ to 1}
\STATE{Set weights of edges in $\edgesactive$ to 2}
\STATE{Set weights of edges in $\edgesdeactivated$ to 3}
\STATE{Find a minimum spanning tree $\tree_{\ell}$ on $\graphl$ (e.g. using Kruskal's algorithm)}
\STATE{$\multileveltree_{L} \leftarrow \multileveltree_{L} \cup \treeactive$}
\ENDFOR
\RETURN $\multileveltree_{L}$
\end{algorithmic}
\end{algorithm}

\begin{remark}
The multi-level tree is not really a tree, because the union of Greville subgrids from different levels does not
form a valid graph. We have decided to call this object a multi-level tree because it gives an idea of how it is 
constructed, and to maintain the \emph{tree-cotree} nomenclature.
\end{remark}

\begin{remark} \label{rem:deactivated_not_needed}
We note that basis functions associated to deactivated edges $\edgeSet_\ell^\deactivated$ are not part of the hierarchical basis $\hbasis^1_{L}$, 
and therefore they are not relevant for the final multi-level tree and cotree. In fact, in our configuration
of a simply connected domain $\Omega \subset \mathbb{R}^2$ with homogeneous boundary conditions, step 3 of the algorithm 
could be removed. \rv{Actually, the algorithm also works in the three-dimensional case, under similar assumptions of a domain with trivial topology.}
However, the presence of deactivated functions is important if the domain has non-trivial topology. 
This will be the subject of a forthcoming paper, which is now in preparation, \rv{and that will also cover the three-dimensional case}.
\end{remark}

\rv{
\begin{remark}
We also note that, for domains with trivial topology and homogeneous boundary conditions, the algorithm is valid independently of the topology of the refined regions. 
That is, even if the refinement creates an artificial hole in the region $G_{\ell,\ell} \setminus G_{\ell,\ell+1}$, the algorithm gives a correct tree-cotree decomposition. We will see an example in the numerical tests of \Cref{sec:numResults}. 
\end{remark}
}

\rv{
\begin{remark}
The correctness of the tree-cotree decomposition for tensor-product splines as in \cite{Kapidani_2022aa} was easily proved, based on the isomorphisms with finite elements \eqref{eq:comm} and the existing theory. Instead, the correctness of the tree-cotree decomposition for the multi-level tree requires a more complex analysis because we are not aware of any result for finite elements with hanging nodes. The theoretical proof of the correctness of the decomposition will be presented in a forthcoming paper, along with the extension to domains with general topology.
\end{remark}
}


\section{Numerical results}
\label{sec:numResults}
To numerically verify the correctness of the proposed tree-cotree decomposition for hierarchical B-splines, we solve
the Maxwell eigenvalue problem in terms of the electric field strength $\Efield$. We write directly the discrete problem
given in weak formulation by
\begin{align*}
 \int_\Omega \nu \curl \, \Efield_h  \cdot \curl \, \Efield'_h = \omega^2 \int_\Omega \epsilon \Efield_h \cdot \Efield'_h, \text{ for all }
 \Efield'_h \in V_h
\end{align*}
where $\Efield_h \in V_h$
and $\omega \in \mathbb{R}$ are the unknown eigenmodes and eigenfrequencies, respectively,
and the material properties $\nu$ and $\epsilon$ represent the magnetic reluctivity and the electric permittivity,
\rv{and for simplicity we take both equal to one}. The discrete space $V_h \subset \mathbf{H}_0(\mathbf{curl};\Omega)$
is chosen as the space of curl-conforming hierarchical splines, $V_h = W^1_L$.

For the solution of this problem with tree-cotree gauging,
we follow the procedure explained in \cite{Manges_1995aa}. While the use of tree-cootree gauging is not
particularly advantageous for the Maxwell eigenvalue problem, we chose it as a model problem because it serves to test the
validity of the gauge. Indeed, when the eigenvalue problem is solved without a gauging method, one obtains zero
eigenvalues associated to gradient fields, and their number corresponds to the dimension of the discrete kernel of the
curl operator. If the problem is gauged correctly, these zero eigenvalues are removed from the solution without affecting
the remaining eigenvalues. If gauging is not applied correctly, spurious modes appear in the computed solution.

We present some numerical tests to test the validity of the tree-cotree decomposition. All the tests are performed using
the GeoPDEs library \cite{Vazquez_2016aa} in MATLAB 2025a, and the eigenvalue problem is solved using the \texttt{eig} command.

\begin{figure}
 \begin{tikzpicture}
	\begin{axis}[width=0.9\linewidth,height=0.4\linewidth,
		xlabel=$i$,
		ylabel=eigenvalue $\lambda_{i}$,
		legend pos=north west]
        \addplot [blue, mark=diamond*, mark options={solid}, only marks] table [x index=0, y index=1, col sep=comma] {figures/data/square_8x8_deg1_eigFull_eigTC_diff.csv};
        \addplot [red, mark=asterisk, only marks] table [x index=0, y index=2, col sep=comma] {figures/data/square_8x8_deg1_eigFull_eigTC_diff.csv};
        \legend{tree-cotree gauged formulation, ungauged formulation}
	\end{axis}
 \end{tikzpicture}
 \caption{Eigenvalues of the Maxwell eigenvalue problem of \Cref{ex:square_8x8} for $p=1$ computed applying tree-cotree gauging and non-zero eigenvalues of the ungauged formulation.}
 \label{fig:square_deg1_eig}
\end{figure}
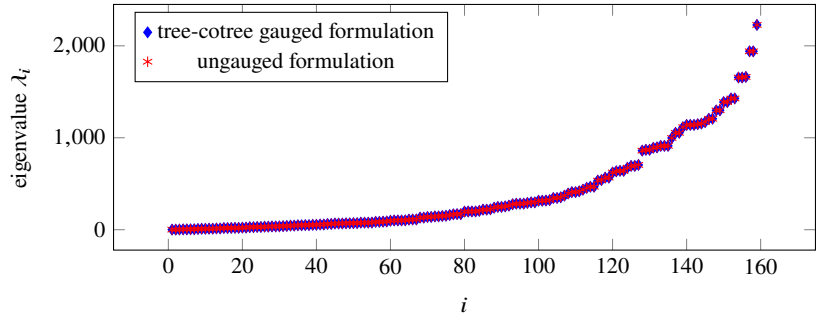

\begin{figure}
 \begin{subfigure}{0.24\textwidth}
  \resizebox{\linewidth}{!}{
  \includegraphics{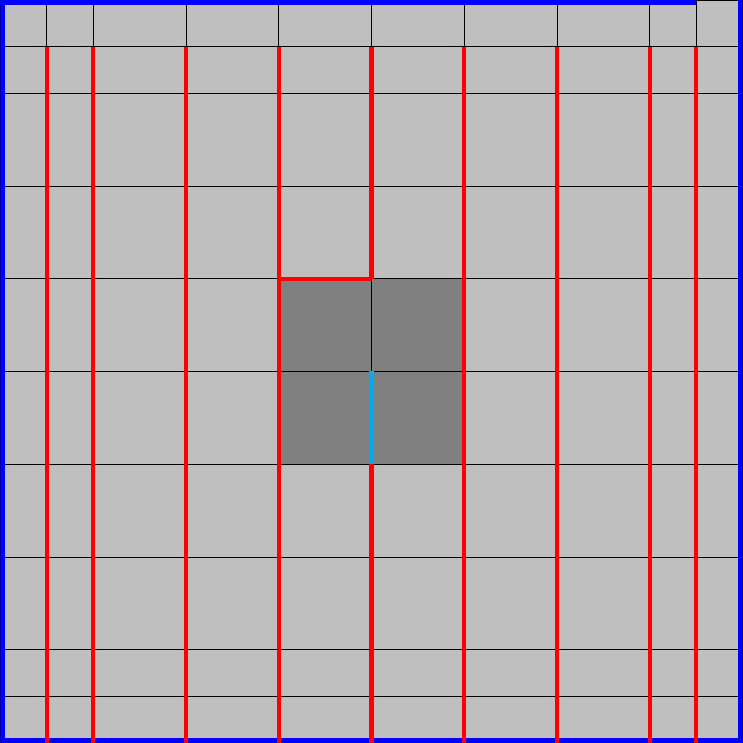}}
  \caption{Tree of level 0}
 \end{subfigure}
 \begin{subfigure}{0.24\textwidth}
  \resizebox{\linewidth}{!}{
  \includegraphics{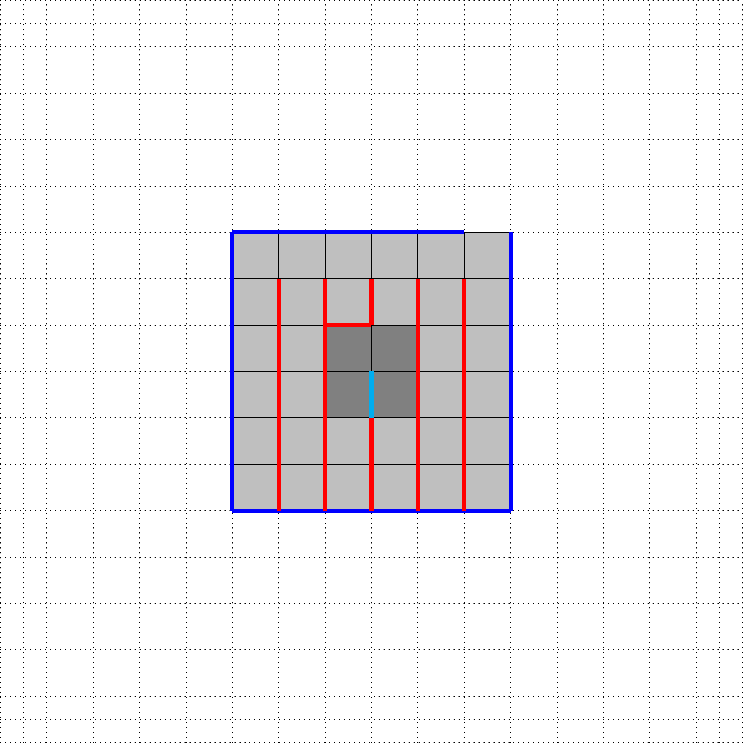}}
  \caption{Tree of level 1}
 \end{subfigure}
 \begin{subfigure}{0.24\textwidth}
  \resizebox{\linewidth}{!}{
  \includegraphics{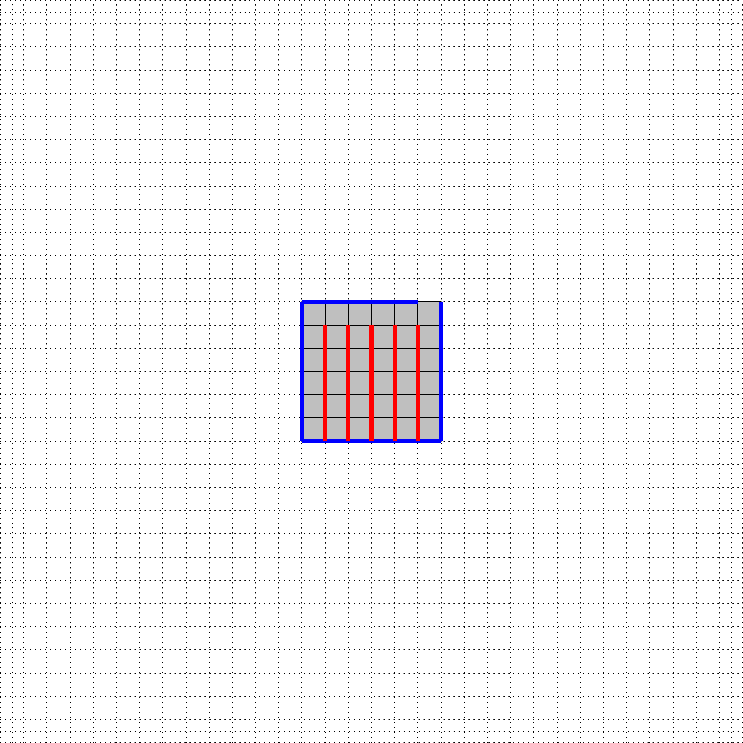}}
  \caption{Tree of level 2}
 \end{subfigure}
 \begin{subfigure}{0.24\textwidth}
  \resizebox{\linewidth}{!}{
  \includegraphics{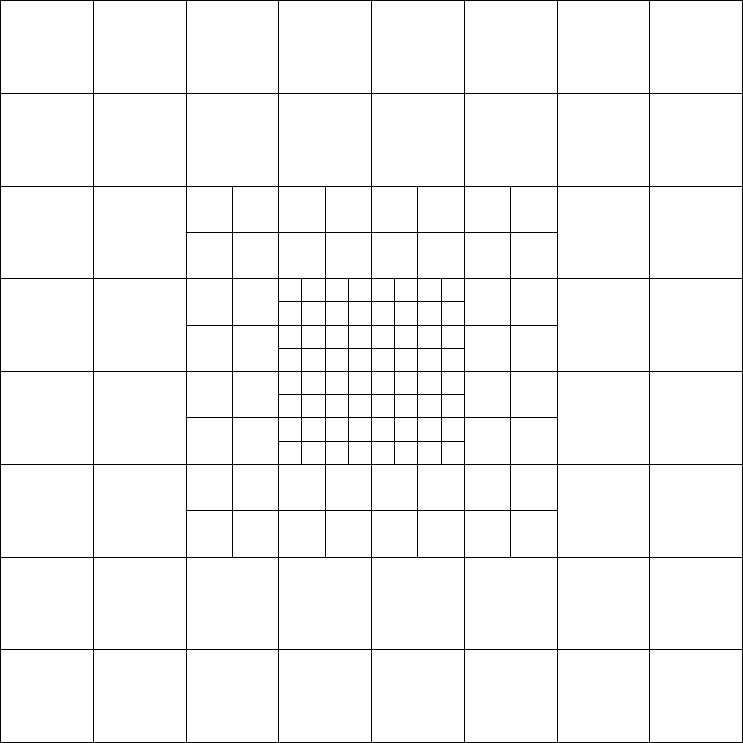}}
  \caption{Hierarchical mesh}
 \end{subfigure}
 \caption{Illustration of the trees built on each level for degree $p=3$, for \Cref{ex:square_8x8}. The first three figures (a)-(c) show
 the tree of each single level, where the whole gray region corresponds to $G_{\ell,\ell}$, and the dark gray region corresponds
 to $G_{\ell,\ell+1}$. The blue edges correspond to the edges on the boundary of $G_{\ell,\ell}$, the red edges to the active
 edges, and the cyan edges to the deactivated edges. The right figure (d) shows the hierarchical mesh.}
 \label{fig:square_deg3}
\end{figure}

\begin{example}[Simple refinement in a square] \label{ex:square_8x8}
    In the first numerical example the domain is the square $\Omega = (0,\pi)^2$, and we
    consider a simple refinement in the central part of the domain. Starting from an $8\times 8$ mesh, we refine at each
    step the central part consisting of a square of $4\times 4$ elements of the finest level, as in \Cref{fig:hier_TC}.
    We stop the refinement at the third level.
    We solve the Maxwell eigenvalue problem with degrees $p=1$ and $p=3$, and compare the eigenvalues computed applying tree-cotree
    gauging with the non-zero eigenvalues obtained with the ungauged formulation, which are considered to be correct.

    The multi-level tree for degree $p=1$ is the one already displayed in \Cref{fig:hier_TC}. The computed eigenvalues, that are shown in \Cref{fig:square_deg1_eig}, demonstrate that the construction of the multi-level tree provides a valid gauge.
    For degree $p=3$, we plot separately the trees for the three levels in \Cref{fig:square_deg3}. As already mentioned, a plot of the
    multi-level tree would not be clear, because the Greville subgrids of the different levels overlap. The comparison of the
    computed eigenvalues in \Cref{fig:square_deg3_eig} confirms the validity of the gauge also for high degree.
\end{example}

\begin{figure}
 \begin{subfigure}{0.24\textwidth}
  \resizebox{\linewidth}{!}{
  \includegraphics{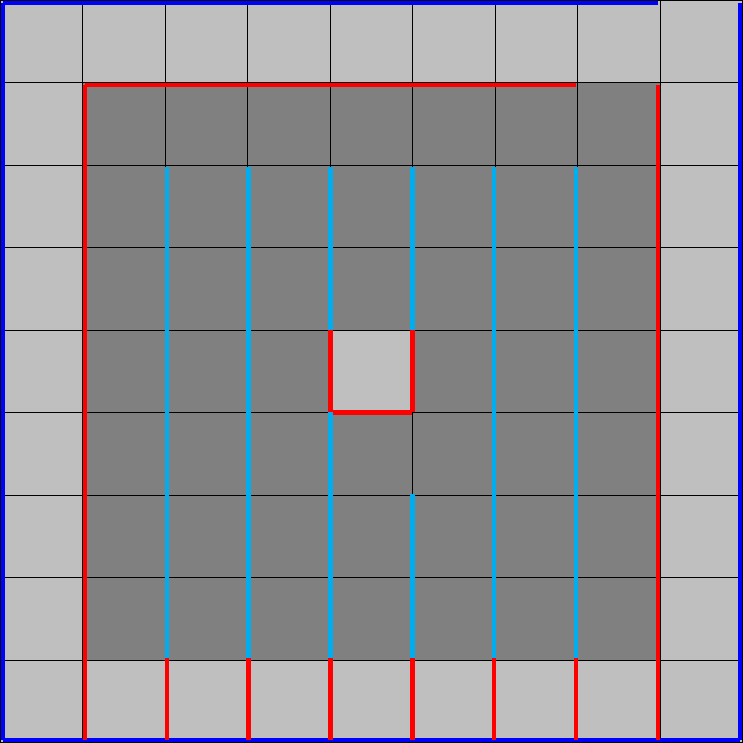}}
  \caption{Tree of level 0}
 \end{subfigure}
 \begin{subfigure}{0.24\textwidth}
  \resizebox{\linewidth}{!}{
  \includegraphics{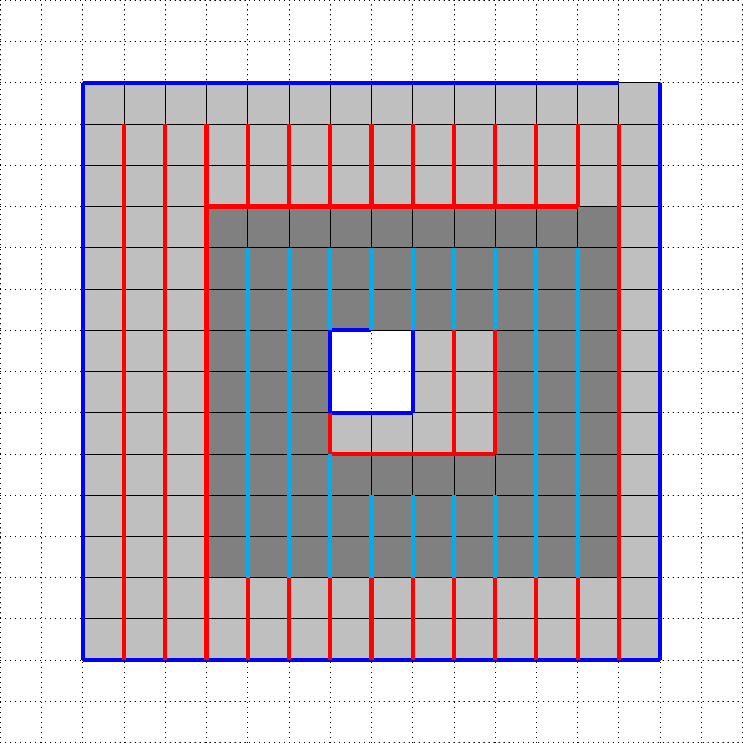}}
  \caption{Tree of level 1}
 \end{subfigure}
 \begin{subfigure}{0.24\textwidth}
  \resizebox{\linewidth}{!}{
  \includegraphics{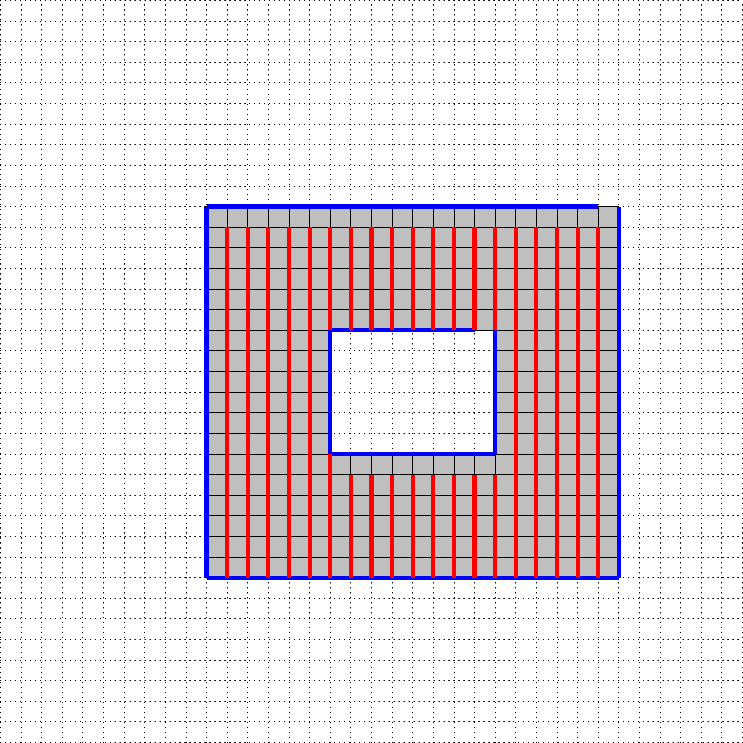}}
  \caption{Tree of level 2}
 \end{subfigure}
  \begin{subfigure}{0.24\textwidth}
  \resizebox{\linewidth}{!}{
  \includegraphics{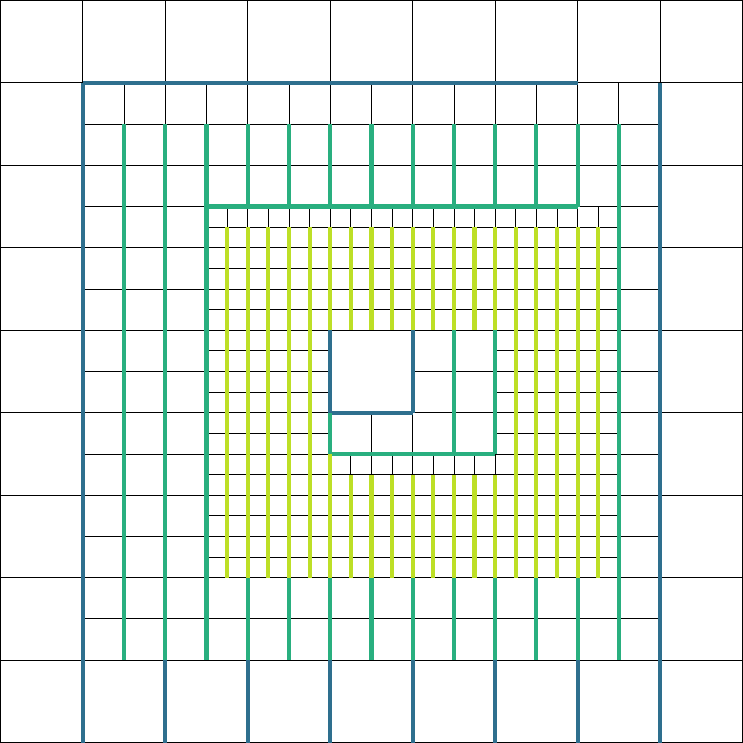}}
  \caption{Multi-level tree}
 \end{subfigure}
 \caption{Illustration of the trees built on each level for degree $p=1$, for \Cref{ex:fake_hole}. The first three figures (a)-(c) show
 the tree of each single level, where the whole gray region corresponds to $G_{\ell,\ell}$, and the dark gray region corresponds
 to $G_{\ell,\ell+1}$. The blue edges correspond to the edges on the boundary of $G_{\ell,\ell}$, the red edges to the active
 edges, and the cyan edges to the deactivated edges. The right figure (d) shows the multi-level tree, combining
 edges of the three levels.}
 \label{fig:fake_hole_deg1}
\end{figure}

\begin{figure}
 \begin{subfigure}{0.24\textwidth}
  \resizebox{\linewidth}{!}{
  \includegraphics{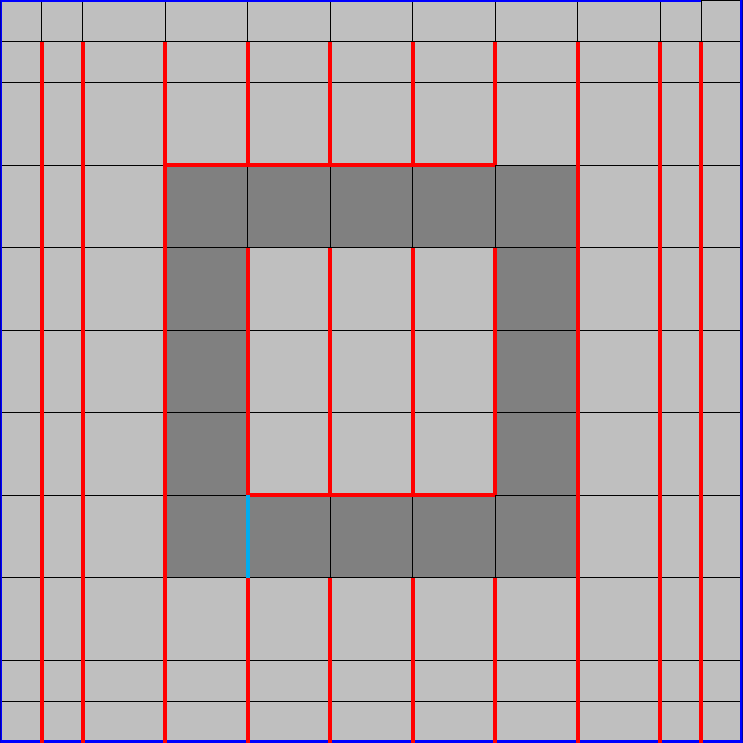}}
  \caption{Tree of level 0}
 \end{subfigure}
 \begin{subfigure}{0.24\textwidth}
  \resizebox{\linewidth}{!}{
  \includegraphics{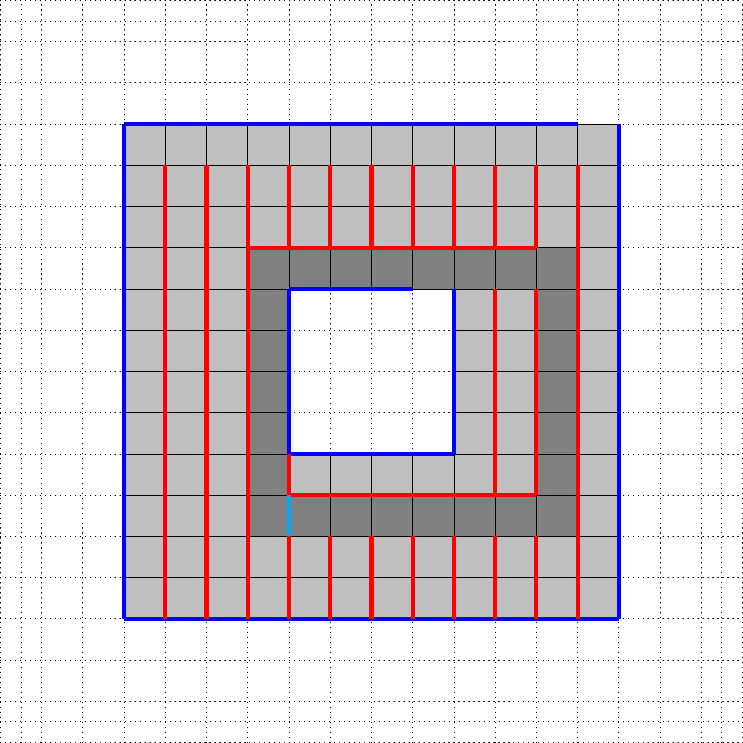}}
  \caption{Tree of level 1}
 \end{subfigure}
 \begin{subfigure}{0.24\textwidth}
  \resizebox{\linewidth}{!}{
  \includegraphics{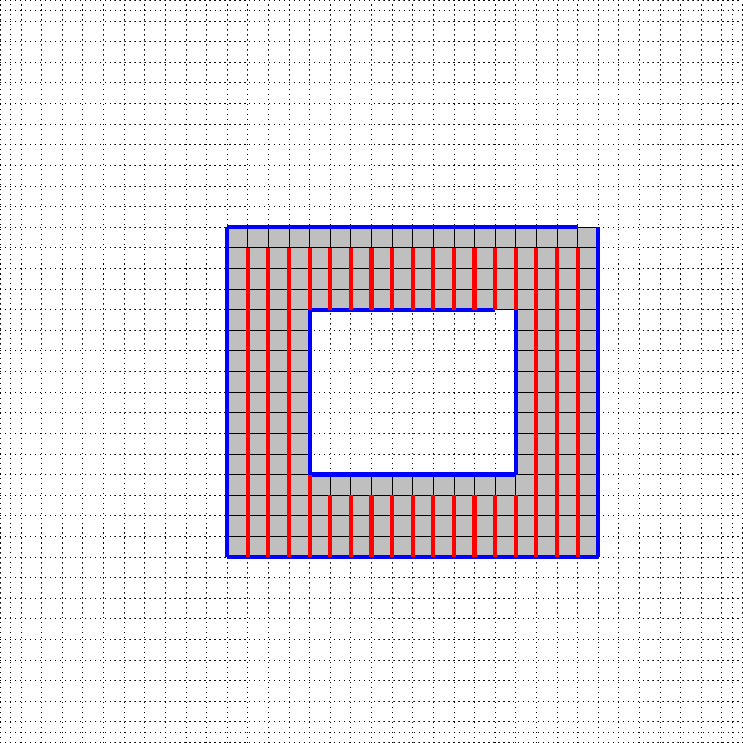}}
  \caption{Tree of level 2}
 \end{subfigure}
 \begin{subfigure}{0.24\textwidth}
  \resizebox{\linewidth}{!}{
  \includegraphics{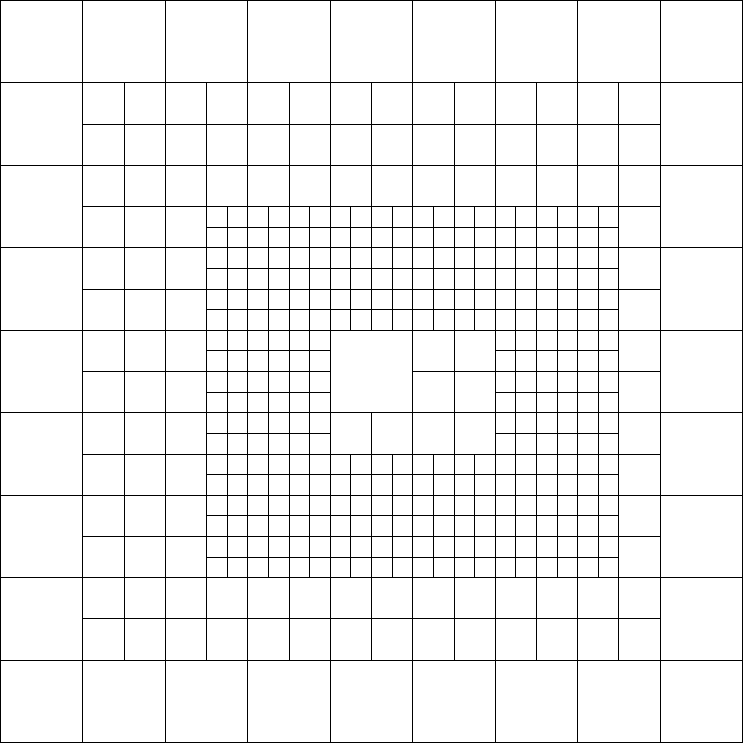}}
  \caption{Hierarchical mesh}\label{fig:fake_hole_hier}
 \end{subfigure}
 \caption{Illustration of the trees built on each level for degree $p=3$, for \Cref{ex:fake_hole}. The first three figures (a)-(c) show
 the tree of each single level, where the whole gray region corresponds to $G_{\ell,\ell}$, and the dark gray region corresponds
 to $G_{\ell,\ell+1}$. The blue edges correspond to the edges on the boundary of $G_{\ell,\ell}$, the red edges to the active
 edges, and the cyan edges to the deactivated edges. The right figure (d) shows the hierarchical mesh.}
 \label{fig:fake_hole_deg3}
\end{figure}

\begin{example}[Refinement creating a hole] \label{ex:fake_hole}
    In the second numerical test we consider again the square $\Omega = (0,\pi)^2$. Starting from an $9 \times 9$ mesh,
    we now refine the elements in such a way that the refined region has a hole in it, see \Cref{fig:fake_hole_hier}, and
    therefore the Greville subgrid of fine levels has non-trivial topology. We refine the mesh in such a way that it only has three
    levels. We plot in \Cref{fig:fake_hole_deg1} and \Cref{fig:fake_hole_deg3} the trees obtained for degrees $p=1$ and $p=3$, respectively. We note that,
    although the Greville subgrids have holes in them, it is not necessary to apply any special treatment due to the topology,
    because all the boundaries are treated as Dirichlet boundary conditions. In fact, the spanning tree connects the two boundaries
    without issues.
    In \Cref{fig:fake_hole_deg1_eig} and \Cref{fig:fake_hole_deg3_eig} we plot the computed eigenvalues for the same degrees. As in the previous tests, there
    is no difference between the eigenvalues obtained for the problem without gauging or applying tree-cotree gauging.
    We can therefore conclude that the tree-cotree decomposition that we propose is correct also in this case.
\end{example}

\begin{figure}
 \begin{subfigure}{0.24\textwidth}
  \resizebox{\linewidth}{!}{
  \includegraphics{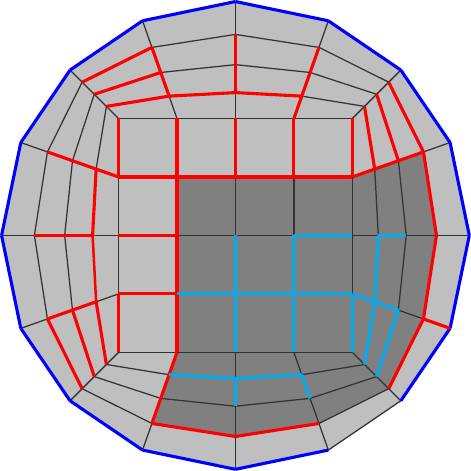}}
  \caption{Tree of level 0}
 \end{subfigure}
 \begin{subfigure}{0.24\textwidth}
  \resizebox{\linewidth}{!}{
  \includegraphics{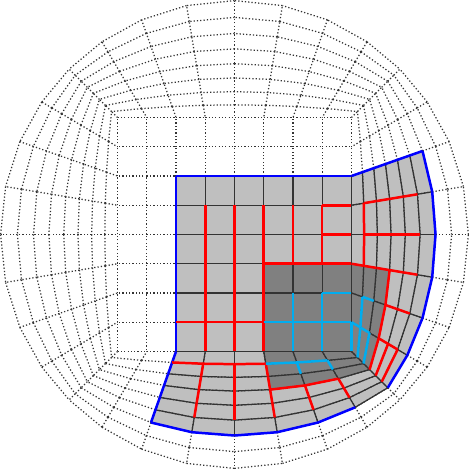}}
  \caption{Tree of level 1}
 \end{subfigure}
 \begin{subfigure}{0.24\textwidth}
  \resizebox{\linewidth}{!}{
  \includegraphics{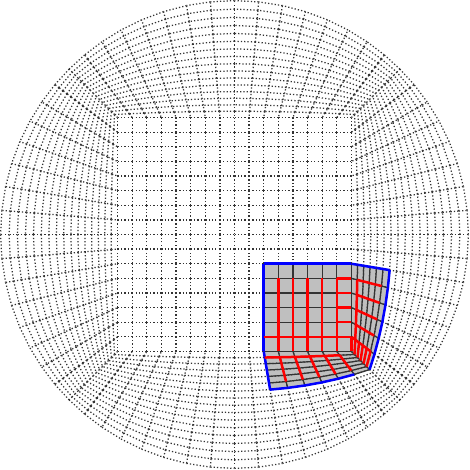}}
  \caption{Tree of level 2}
 \end{subfigure}
  \begin{subfigure}{0.24\textwidth}
  \resizebox{\linewidth}{!}{
  \includegraphics{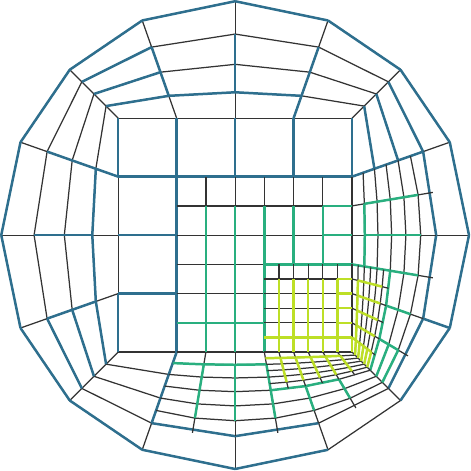}}
  \caption{Multi-level tree}
 \end{subfigure}
 \caption{\mm{Illustration of the trees built on each level for degree $p=1$, for \Cref{ex:circle}. The first three figures (a)-(c) show
 the tree of each single level, where the whole gray region corresponds to $G_{\ell,\ell}$, and the dark gray region corresponds
 to $G_{\ell,\ell+1}$. The blue edges correspond to the edges on the boundary of $G_{\ell,\ell}$, the red edges to the active
 edges, and the cyan edges to the deactivated edges. The right figure (d) shows the multi-level tree, combining
 edges of the three levels.}}
 \label{fig:circle_deg1}
\end{figure}

\begin{figure}
 \begin{subfigure}{0.24\textwidth}
  \resizebox{\linewidth}{!}{
  \includegraphics{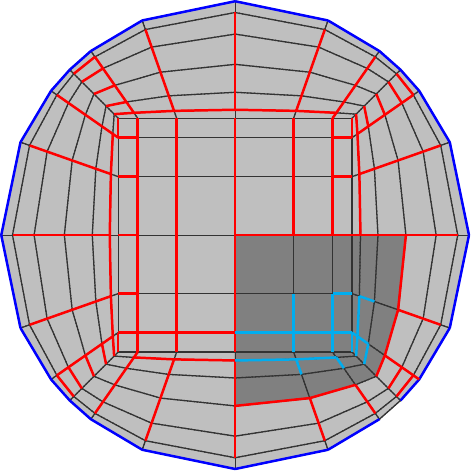}}
  \caption{Tree of level 0}
 \end{subfigure}
 \begin{subfigure}{0.24\textwidth}
  \resizebox{\linewidth}{!}{
  \includegraphics{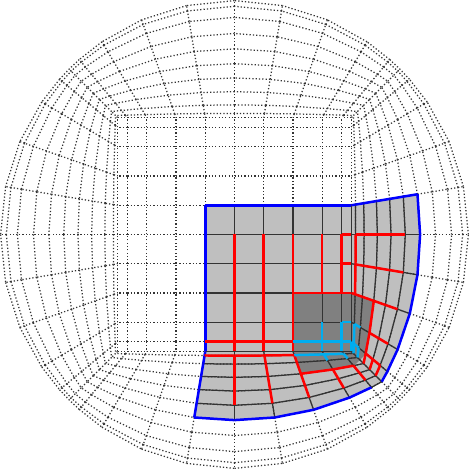}}
  \caption{Tree of level 1}
 \end{subfigure}
 \begin{subfigure}{0.24\textwidth}
  \resizebox{\linewidth}{!}{
  \includegraphics{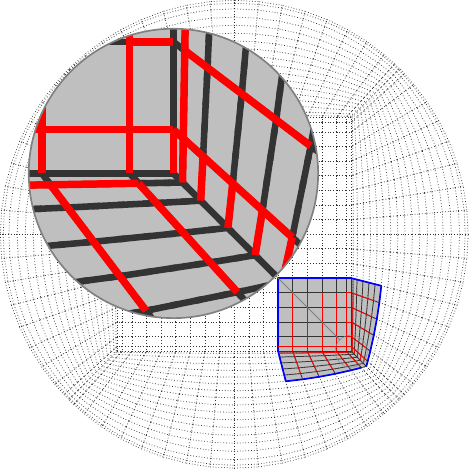}}
  \caption{Tree of level 2}
 \end{subfigure}
  \begin{subfigure}{0.24\textwidth}
  \resizebox{\linewidth}{!}{
  \includegraphics{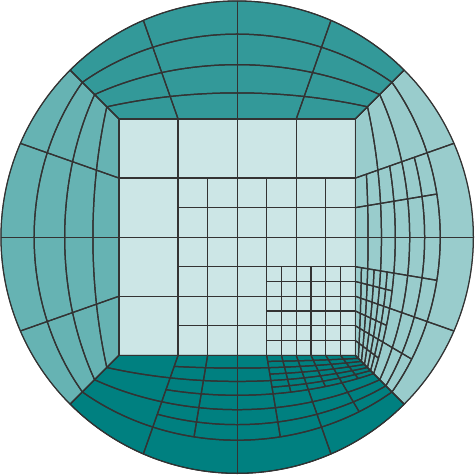}}
  \caption{Hierarchical mesh}\label{fig:circle_hier}
 \end{subfigure}
 \caption{\mm{Illustration of the trees built on each level for degree $p=3$, for \Cref{ex:circle}. The first three figures (a)-(c) show
 the tree of each single level, where the whole gray region corresponds to $G_{\ell,\ell}$, and the dark gray region corresponds
 to $G_{\ell,\ell+1}$. The blue edges correspond to the edges on the boundary of $G_{\ell,\ell}$, the red edges to the active
 edges, and the cyan edges to the deactivated edges. The right figure (d) shows the hierarchical mesh and the five patches are shown in different colors.}}
 \label{fig:circle_deg3}
\end{figure}

\begin{example}[Multi-patch circle] \label{ex:circle}
\mm{
    In the third numerical test the domain is a disk of radius $R=1$. It consists of five patches: a square patch $\Omega^{(1)} = (0,1)^2$ in the center and four patches $\Omega^{(2)},\Omega^{(3)},\Omega^{(4)},\Omega^{(5)}$ around it.
    We consider a refinement with three levels around one of the nodes where three patches meet. The patches and refinement are shown in \Cref{fig:circle_hier}.
    We plot in \Cref{fig:circle_deg1} and \Cref{fig:circle_deg3} the trees obtained for degrees $p=1$ and $p=3$, respectively.
    In \Cref{fig:circle_deg1_eig} and \Cref{fig:circle_deg3_eig} we plot the computed eigenvalues for the same degrees. As in the previous tests, there
    is no difference between the eigenvalues obtained for the problem without gauging or applying tree-cotree gauging.
    We can therefore conclude that the tree-cotree decomposition that we propose is also correct for the multi-patch case.
    }
\end{example}

%

\begin{figure}
 \begin{tikzpicture}
	\begin{axis}[width=0.9\linewidth,height=0.4\linewidth,
		xlabel=$i$,
		ylabel=eigenvalue $\lambda_{i}$,
		legend pos=north west]
        \addplot [blue, mark=diamond*, mark options={solid}, only marks] table [x index=0, y index=1, col sep=comma] {figures/data/square_8x8_deg3_eigFull_eigTC_diff.csv};
        \addplot [red, mark=asterisk, only marks] table [x index=0, y index=2, col sep=comma] {figures/data/square_8x8_deg3_eigFull_eigTC_diff.csv};
        \legend{tree-cotree gauged formulation, ungauged formulation}
	\end{axis}
 \end{tikzpicture}
 \caption{Eigenvalues of the Maxwell eigenvalue problem of \Cref{ex:square_8x8} for $p=3$ computed applying tree-cotree gauging and non-zero eigenvalues of the ungauged formulation.}
 \label{fig:square_deg3_eig}
\end{figure}

\begin{figure}
 \begin{tikzpicture}
	\begin{axis}[width=0.9\linewidth,height=0.4\linewidth,
		xlabel=$i$,
		ylabel=eigenvalue $\lambda_{i}$,
		legend pos=north west]
        \addplot [blue, mark=diamond*, mark options={solid}, only marks] table [x index=0, y index=1, col sep=comma] {figures/data/fake_hole_three_levels_9x9_deg1_eigFull_eigTC_diff.csv};
        \addplot [red, mark=asterisk, only marks] table [x index=0, y index=2, col sep=comma] {figures/data/fake_hole_three_levels_9x9_deg1_eigFull_eigTC_diff.csv};
        \legend{tree-cotree gauged formulation, ungauged formulation}
	\end{axis}
 \end{tikzpicture}
 \caption{Eigenvalues of the Maxwell eigenvalue problem of \Cref{ex:fake_hole} for $p=1$ computed applying tree-cotree gauging and non-zero eigenvalues of the ungauged formulation.}
 \label{fig:fake_hole_deg1_eig}
\end{figure}

\begin{figure}
 \begin{tikzpicture}
	\begin{axis}[width=0.9\linewidth,height=0.4\linewidth,
		xlabel=$i$,
		ylabel=eigenvalue $\lambda_{i}$,
		legend pos=north west]
        \addplot [blue, mark=diamond*, mark options={solid}, only marks] table [x index=0, y index=1, col sep=comma] {figures/data/fake_hole_three_levels_9x9_deg3_eigFull_eigTC_diff.csv};
        \addplot [red, mark=asterisk, only marks] table [x index=0, y index=2, col sep=comma] {figures/data/fake_hole_three_levels_9x9_deg3_eigFull_eigTC_diff.csv};
        \legend{tree-cotree gauged formulation, ungauged formulation}
	\end{axis}
 \end{tikzpicture}
 \caption{Eigenvalues of the Maxwell eigenvalue problem of \Cref{ex:fake_hole} for $p=3$ computed applying tree-cotree gauging and non-zero eigenvalues of the ungauged formulation.}
 \label{fig:fake_hole_deg3_eig}
\end{figure}

\begin{figure}
 \begin{tikzpicture}
	\begin{axis}[width=0.9\linewidth,height=0.4\linewidth,
		xlabel=$i$,
		ylabel=eigenvalue $\lambda_{i}$,
		legend pos=north west]
        \addplot [blue, mark=diamond*, mark options={solid}, only marks] table [x index=0, y index=1, col sep=comma] {figures/data/circle_mp_three_levels_deg1_eigFull_eigTC_diff.csv};
        \addplot [red, mark=asterisk, only marks] table [x index=0, y index=2, col sep=comma] {figures/data/circle_mp_three_levels_deg1_eigFull_eigTC_diff.csv};
        \legend{tree-cotree gauged formulation, ungauged formulation}
	\end{axis}
 \end{tikzpicture}
 \caption{\mm{Eigenvalues of the Maxwell eigenvalue problem of \Cref{ex:circle} for $p=1$ computed applying tree-cotree gauging and non-zero eigenvalues of the ungauged formulation.}}
 \label{fig:circle_deg1_eig}
\end{figure}

\begin{figure}
 \begin{tikzpicture}
	\begin{axis}[width=0.9\linewidth,height=0.4\linewidth,
		xlabel=$i$,
		ylabel=eigenvalue $\lambda_{i}$,
		legend pos=north west]
        \addplot [blue, mark=diamond*, mark options={solid}, only marks] table [x index=0, y index=1, col sep=comma] {figures/data/circle_mp_three_levels_deg3_eigFull_eigTC_diff.csv};
        \addplot [red, mark=asterisk, only marks] table [x index=0, y index=2, col sep=comma] {figures/data/circle_mp_three_levels_deg3_eigFull_eigTC_diff.csv};
        \legend{tree-cotree gauged formulation, ungauged formulation}
	\end{axis}
 \end{tikzpicture}
 \caption{\mm{Eigenvalues of the Maxwell eigenvalue problem of \Cref{ex:circle} for $p=3$ computed applying tree-cotree gauging and non-zero eigenvalues of the ungauged formulation.}}
 \label{fig:circle_deg3_eig}
\end{figure}

%
%
\section*{Acknowledgements}
The work of Rafael Vázquez has been partially funded by 
Consellería de Educación, Ciencia, Universidades e Formación Profesional - Xunta de Galicia (2025-AD059 and ED431F 2025/03),
and by the Spanish State Research Agency (PID2024-156071NB-I00).
The work of Melina Merkel is supported
by the joint DFG/FWF Collaborative Research Centre CREATOR (DFG: Project-ID 492661287/TRR 361; FWF: 10.55776/F90) at TU Darmstadt, TU Graz and JKU Linz,
by the LOEWE Project 1450/23-04 via the Hessian Ministry of Science and Art and the Hessen-Agentur
and by
the Association of Friends of TU Darmstadt via an Ernst Ludwig Mobility Grant.
\bibliographystyle{spmpsci}
\bibliography{bibtex,references}

\end{document}